\documentclass[a4paper,11pt]{article}

\usepackage{bm}
\usepackage{anysize}
\usepackage{amsthm}
\usepackage{amssymb}
\usepackage{amsmath}
\usepackage{fancybox}
\usepackage{amsfonts}
\usepackage[T1]{fontenc}
\usepackage{latexsym}
\usepackage{makeidx}
\usepackage{cite}
\usepackage[numbers]{natbib}

\usepackage[pdftex]{graphicx}
\usepackage[pdftex]{color}

\usepackage[all]{xy}

\newtheorem{teo}{Theorem}

\newtheorem{pro}{Proposition}
\newtheorem{cor}{Corollary}
\newtheorem{defi}{Definition} 

\newcommand{\fd}{\rightarrow}

\newcommand{\inc}{\subset}

\newcommand{\ba}{\overline}
\newcommand{\onda}{\widetilde}

\newcommand{\al}{\alpha}
\newcommand{\be}{\beta}
\newcommand{\lan}{\lambda}
\newcommand{\fhi}{\varphi}
\newcommand{\del}{\delta}
\newcommand{\Del}{\Delta}
\newcommand{\gam}{\gamma}
\newcommand{\Gam}{\Gamma}

\newcommand{\Om}{\Omega}

\newcommand{\Z}{\mathbb{Z}}
\newcommand{\N}{\mathbb{N}}

\newcommand{\R}{\mathbb{R}}
\newcommand{\C}{\mathbb{C}}

\newcommand{\Sa}{\mathbb{S}}

\newcommand{\pa}{\partial}

\newcommand{\p}{\grave{}}
\newcommand{\f}{\mathfrak{f}}

\newtheorem{remark}{Remark}[section]

\def\pf{\par\noindent {\em Proof.}~\par\noindent}

\def\lim{\mathop{\mbox{\normalfont lim}}\limits}

\def\pf{\par\noindent {\em Proof. }}

\def\pa{\partial}

\def\ve{\varepsilon}

\begin{document}

\date{}

\title{Bochner-Martinelli formula in superspace}
\small{
\author
{Juan Bory Reyes$^1$, Al\'i Guzm\'an Ad\'an$^2$, Frank Sommen$^2$}
\vskip 1truecm
\date{\small  $^1${SEPI-ESIME-ZAC-Instituto Polit\'ecnico Nacional, Ciudad M\'exico 07738, M\'exico}\\
 $^2${Clifford Research Group, Department of Mathematical Analysis, Faculty of Engineering and Architecture, Ghent University, Krijgslaan 281, 9000 Gent, Belgium. }}}

\maketitle

\begin{abstract}

In a series of recent papers, a harmonic and hypercomplex function theory in superspace has been established and amply developed. In this paper, we address the problem of establishing Cauchy integral formulae in the framework of Hermitian Clifford analysis in superspace. This allows us to obtain a successful extension of the classical Bochner-Martinelli formula to superspace by means of the corresponding projections on the space of spinor-valued superfunctions.

\noindent

\vspace{0.3cm}

\small{ }

\noindent
\textbf{}
\end{abstract}

\section{Introduction}

The Bochner-Martinelli integral representation constitutes a classical generalization, to the case of several complex variables, of the Cauchy integral formula for holomorphic functions in the complex plane. This representation reads for every holomorphic complex function $f$ on some bounded domain $\Om\inc \C^m$, with smooth boundary $\pa\Om$, as
\begin{align}\label{B-Moriginal}
f(\underline{U}) &= \int_{\pa\Om} f(\underline{Z}) \;  \mathcal K(\underline{Z}, \underline{U}),  & \underline{U}\in \Om,
\end{align}
where $\mathcal K(\underline{Z}, \underline{U})$ is the exterior differential form of type ($m,m-1$) given by
\[
\mathcal K(\underline{Z}, \underline{U}) = \frac{(m-1)!}{(2\pi i)^m} \; \sum_{j=1}^m (-1)^{j-1}\frac{z^c_j-u^c_j}{|\underline{Z}-\underline{U}|^{2m}} \; dz_1^c\wedge \cdots \wedge dz_{j-1}^c\wedge dz_{j+1}^c \wedge \cdots \wedge dz_m^c \wedge dz_1\wedge \cdots \wedge dz_m,
\]
and $\cdot^c$ denotes the complex conjugation. The form $\mathcal K(\underline{Z}, \underline{U})$ is the so-called Bochner-Martinelli kernel. When $m=1$, this kernel reduces to the Cauchy kernel $(2\pi i)^{-1}(z-u)^{-1}\, dz$, {whence} formula (\ref{B-Moriginal}) reduces to the traditional Cauchy integral formula in one complex variable. For $m>1$,  $\mathcal K(\underline{Z}, \underline{U})$ fails to be holomorphic but it still remains harmonic, see e.g.\ \cite{MR1409816}. Formula (\ref{B-Moriginal}) was obtained independently and through different methods by Martinelli and Bochner, see e.g.\ \cite{MR1162310} for a detailed description. The interest for proving different generalizations of the classical Bochner-Martinelli formula has emerged as a successful research topic.

A second important generalization of the Cauchy integral formula is offered by so-called Euclidean Clifford analysis, where functions defined in Euclidean space $\R^{2m}$ and taking values in a Clifford algebra are considered. This theory is focussed on the notion of monogenic functions, i.e. null solutions of the Dirac operator $\pa_{\underline x}$ factorizing the Laplace operator: $\pa_{\underline x}^2=-\Del_{2m}$. Standard references on this setting are \cite{MR697564, MR1169463, MR1130821}. In this framework the Clifford-Cauchy integral formula for monogenic functions reads
\begin{align*}
f(\underline{y}) &= \int_{\pa \Om} \fhi_1^{2m|0} (\underline{x} - \underline{y}) \, n(\underline{x}) \, f(\underline{x})\, dS_{\underline{x}},  &  \underline{y}&\in \Om,
\end{align*}
where $ \fhi_1^{2m|0} (\underline{x} - \underline{y}) = \frac{-1}{|\Sa^{2m-1}|} \frac{\underline{x} - \underline{y}}{|\underline{x} - \underline{y}|^{2m}}$ is the so-called Cauchy kernel, $|\Sa^{2m-1}|$ is the area of the unit sphere $\Sa^{2m-1}$ in $\R^{2m}$,  $n(\underline{x})$ denotes the exterior normal vector to $\pa \Om$ at the point $\underline{x}\in \pa \Om$ and $dS_{\underline{x}}$ is the Lebesgue surface measure in $\pa \Om$. This formula {has been} a cornerstone in the development of the monogenic function theory. 

{Both integral representations above} were proven to be related when one considers so-called Hermitian Clifford analysis, which constitutes yet a refinement of the Euclidean case. This refinement focusses on the simultaneous null solutions of the complex Hermitian Dirac operators $\pa_{\underline{Z}}$ and $\pa_{\underline{Z}^\dag}$ which decomposes the Laplace operator in the sense that $4(\pa_{\underline{Z}}\pa_{\underline{Z}^\dag}+ \pa_{\underline{Z}^\dag} \pa_{\underline{Z}}) =\Del_{2m}$. We refer the reader to \cite{MR2490567, Brackx2007, MR2465736, MMA:MMA378} for a general overview. {Indeed, in \cite{MR2540516}} a Cauchy integral formula for Hermitian monogenic functions was obtained by passing to the framework of circulant ($2\times 2$) matrix functions. {This} Hermitian Cauchy integral representation was proven to reduce to the traditional Bochner-Martinelli formula (\ref{B-Moriginal}) when considering the special case of functions taking values in the zero-homogeneous part of complex spinor space. This means that the theory of
Hermitian monogenic functions not only refines Euclidean Clifford analysis
(and thus harmonic analysis as well), but also has strong connections with the
theory of functions of several complex variables, even encompassing some of its results. 

Our main goal in the current paper is to extend the Bochner-Martinelli formula (\ref{B-Moriginal}) to superspace by exploiting the above-described relation with Clifford analysis. Superspaces play an important rôle in contemporary theoretical physics, e.g. in the particle theory of supersymmetry, supergravity or superstring theories, etc. Traditionally they have been studied using algebraic and geometrical methods (see e.g.\ \cite{Berezin:1987:ISA:38130, MR565567, MR778559}). More recently harmonic and Clifford analysis have been extended to superspace by introducing some important differential operators (such as Dirac and Laplace operators) and by studying special functions and orthogonal polynomials related to these operators, see e.g.\ \cite{Bie2007, de2007clifford, MR2344451, 1751-8121-42-24-245204, MR2521367}. The basics of Hermitian Clifford analysis in superspace were introduced in \cite{DS_Guz_Somm2} following the notion of an abstract complex structure in the Hermitian radial algebra, developed in \cite{MMA:MMA378, MR3656057}. 
Some particular aspects related to the invariance properties with respect to underlying Lie groups and Lie algebras in this setting has been already studied, see e.g\  \cite{MR3702693}.

 In the current paper, we further develop this function theory by establishing a Cauchy integral representation for Hermitian monogenic superfunctions. To this end, we use a general distributional approach to integration in superspace provided in the works \cite{Guz_Somm5, MR2539324}. This way of integrating constitutes a suitable extension to superspace of H\"ormander's formula which links classical real integration with the use of delta distributions,  see \cite[Theorem 6.1.5]{MR1996773}. As has been shown in the literature, this approach features an important number of advantages. First of all, it contains and unifies of all the previously known approaches, i.e.\ compositions of real integration with the Berezin integral (see e.g\ \cite{Berezin:1987:ISA:38130, MR2521367}), and integration over the supersphere (see e.g.\ \cite{MR2539324, MR2344451}). Secondly, it considers integrals as functionals depending on a fixed {\it phase function}  instead of depending on geometrical sets of points. This leads to a calculus independent of representations of the corresponding bosonic and fermionic variables as co-ordinates defined in a particular superspace over an underlying commutative superalgebra. 
 
 Finally, we establish the connection between Hermitian monogenicity and holomorphicity in superspace by considering a specific class of spinor-valued superfunctions (Section \ref{HolSSSection}). As one may have expected, the obtained (super) Hermitian Cauchy integral formula (\ref{HCauchyFormHerSupAnaMatSpeCa}) reduces, when considering the correct projections, to a new extension of the Bochner-Martinelli formula for holomorphic functions in superspace. 
  
\section{Preliminaries}

Superanalysis or analysis on superspace considers not only commuting (bosonic) but also anti-commuting (fermionic) variables. In this paper we follow the extension of harmonic and Clifford analysis to superspace (see \cite{de2007clifford, Bie2007, MR2344451, 1751-8121-42-24-245204}). This approach considers $m$ commuting variables $x_1,\ldots, x_m$ and $2n$ anti-commuting variables $x\p_1, \ldots, x\p_{2n}$ in a purely symbolic way, i.e.\
\begin{align*}
 x_jx_k &= x_kx_j, & x\p_j x\p_k &= -x\p_k x\p_j, & x_jx\p_k &= x\p_k x_j. 
\end{align*}
This means that $x_1,\ldots, x_m$ are interpreted as generators of the polynomial algebra $\R[x_1,  \ldots, x_m]$ while $x\p_1, \ldots, x\p_{2n}$ generate a Grassmann algebra $\mathfrak{G}_{2n}$. We will denote by $\mathfrak{G}^{(ev)}_{2n}$ and $\mathfrak{G}^{(odd)}_{2n}$ the subalgebras of even and odd elements of $\mathfrak{G}_{2n}$ respectively. The algebra of super-polynomials, i.e.\ polynomials in the variables $x_1,  \ldots, x_m, x\p_1, \ldots, x\p_{2n}$, is defined by
\[ \mathcal{P}:= \mbox{Alg}_{\R}(x_1,  \ldots, x_m, x\p_1, \ldots, x\p_{2n})= \R[x_1,  \ldots, x_m] \otimes \mathfrak{G}_{2n}.\]
The bosonic and fermionic partial derivatives $\pa_{x_j}$, $\pa_{x\p_j}$ are defined as endomorphisms on $\mathcal P$ by the  relations 
\begin{equation*}
\begin{cases} \pa_{x_j}[1]=0,\\
 \pa_{x_j} x_k- x_k \pa_{x_j}=\del_{j,k},\\
 \pa_{x_j} x\p_k=x\p_k \pa_{x_j}, \;\; 
 \end{cases}
\hspace{.5cm} 
\begin{cases} \pa_{x\p_j}[1]=0,\\
\pa_{x\p_j} x\p_k+ x\p_k \pa_{x\p_j}=\del_{j,k}\\
\pa_{x\p_j} x_k=x_k\pa_{x\p_j}, 
\end{cases}
\end{equation*} 
that can be applied recursively. From this definition it immediately follows that $\pa_{x_j} \pa_{x_k}=\pa_{x_k}\pa_{x_j}$, $ \pa_{x\p_j}\pa_{x\p_k}=-\pa_{x\p_k}\pa_{x\p_j}$ and $\pa_{x_j}\pa_{x\p_k}=\pa_{x\p_k}\pa_{x_j}$. 

The flat supermanifold corresponding to the variables $x_j$, $x\p_j$ is denoted by $\R^{m|2n}$. The full algebra of functions on this supermanifold is $C^\infty(\R^m)\otimes \mathfrak{G}_{2n}$ where $C^\infty(\R^m)$ denotes the space of infinitely many differentiable complex-valued functions defined in $\R^m$. The partial derivatives $\pa_{x_j}$, $\pa_{x\p_j}$ extend from $\mathcal P$ to $C^\infty(\R^m)\otimes \mathfrak{G}_{2n}$ by density. 

Clifford algebras in superspace are introduced by $m$ orthogonal Clifford generators  $e_1,\ldots, e_m$ and $2n$ symplectic Clifford generators $e\p_1, \ldots, e\p_{2n}$ which are subjected to the multiplication relations
\begin{align}\label{CommRules}
e_je_k+e_ke_j=-2\del_{j,k}, \hspace{.3cm} e_je\p_k+e\p_ke_j=0, \hspace{.3cm} e\p_je\p_k-e\p_ke\p_j=g_{j,k},
\end{align}
where $g_{j,k}$ is a symplectic form defined by
\[g_{2j,2k}=g_{2j-1,2k-1}=0, \hspace{.5cm} g_{2j-1,2k}=-g_{2k,2j-1}=\del_{j,k}, \hspace{.5cm} j,k=1,\ldots,n.\]
These generators are combined with the algebra of super-polynomials $\mathcal P$  giving rise to the algebra of Clifford valued super-polynomials $\mathcal P \otimes\mathcal C_{m,2n}$; 
where elements in $\mathcal C_{m,2n}:=\mbox{Alg}_\R(e_1,\ldots,e_m,e\p_1,\ldots,e\p_{2n})$ commute  with elements in $\mathcal P$. Also the partial derivatives $\pa_{x_j}$, $\pa_{x\p_j}$ commute with the elements in the algebra $\mathcal C_{m,2n}$. 

 The most important element of the algebra  $\mathcal P \otimes\mathcal C_{m,2n}$ is the supervector variable
\begin{equation}\label{SupVec}
{\bf x}=\underline{x}+\underline{x\p}=\sum_{j=1}^mx_je_j+\sum_{j=1}^{2n}x\p_j e\p_j.
\end{equation}

Functions in $C^\infty(\R^m)\otimes \mathfrak{G}_{2n}$ can be explicitly written as 
\begin{equation}\label{SupFunc}
F({\bf x})=F(\underline{x},\underline{x\p})=\sum_{A\inc \{1, \ldots, 2n\}} \, F_A(\underline{x})\, \underline{x}\p_A, 
\end{equation}
where $F_A(\underline{x})\in C^\infty(\R^m)$ and  $ \underline{x}\p_A$ is defined as $x\p_{j_1}\ldots x\p_{j_k}$ with $A=\{j_1,\ldots, j_k\}$, $1\leq j_1< \ldots< j_k\leq 2n$.  
Similarly, one may consider other spaces of superfunctions of the form $\mathcal F\otimes \mathfrak{G}_{2n}$  where $\mathcal{F}=C^k(\Om), L_2(\Om), \ldots, $ with $\Om\inc \R^m$. In general, the bosonic functions $F_A$ in (\ref{SupFunc}) are complex-valued. We say that $F$ is a {\it real superfunction} when all the elements $F_A$ are real-valued. 

Every superfunction can be written as the sum $F({\bf x})= F_0(\underline{x})+ {\bf F}(\underline x, \underline{x}\p)$ where the  complex-valued function $F_0(\underline{x})=F_\emptyset(\underline{x})$ is called the {\it body} $F$, and ${\bf F}=\sum_{|A|\geq 1} \, F_A(\underline{x})\,\underline{x}\p_A$ is the {\it nilpotent part} of $F$. Indeed, it is clearly seen that  ${\bf F}^{2n+1}=0$.

Following the classical approach, the (real) support $supp\; F$ of a superfunction $F\in C^\infty(\R^m)\otimes \mathfrak{G}_{2n}$ is defined as the closure of the set of all points in $\R^m$ for which the function $F(\cdot, \underline{x}\p):\R^m\fd \mathfrak{G}_{2n}$ is not zero, see (\ref{SupFunc}). From this definition, it immediately follows  that $supp\; F =\bigcup_{A\inc \{1, \ldots, 2n\}} supp\; F_A$.

The  bosonic and fermionic Dirac operators are defined 
by
\[\pa_{\underline x}=\sum_{j=1}^m e_j\pa_{x_j}, \hspace{1cm} \pa_{\underline x\p}=2\sum_{j=1}^n \left(e\p_{2j}\pa_{x\p_{2j-1}}-e\p_{2j-1}\pa_{x\p_{2j}}\right),\]
which lead to the left and right super Dirac operators (super-gradient)
$\pa_{\bf x} \cdot =\pa_{\underline x\p}\cdot -\pa_{\underline x}\cdot$ and $\cdot\pa_{\bf x}  =- \cdot\pa_{\underline x\p} - \cdot \pa_{\underline x}$ respectively. As in the classical Clifford setting, the action of $\pa_{\bf x}$ on the vector variable ${\bf x}$ results in the superdimension $M=m-2n$. {In this paper we work with general superdimension $M$ (but $m\neq 0$).

Given an open set $\Om\inc \R^m$, a function $F\in C^1(\Om)\otimes \mathfrak{G}_{2n}\otimes \mathcal C_{m,2n}$ is said to be (left) {\it super monogenic} if $\pa_{\bf x}[F]=0$. As the super Dirac operator factorizes the super Laplace operator:
\[\Del_{m|2n}=-\pa_{\bf x}^2=\sum_{j=1}^{m} \pa^2_{x_j}-4\sum_{j=1}^n \pa_{x\p_{2j-1}}\pa_{x\p_{2j}},\]
monogenicity also constitutes  a refinement of harmonicity in superanalysis. More details on the theory of super-monogenic and super-harmonic functions can be found for instance in \cite{MR2344451, 1751-8121-42-24-245204, MR2386499, MR2405887, MR3375856}.

 It is possible to produce interesting even superfunctions out of known functions from real analysis. Indeed, consider an analytic function $F\in C^\infty(\R)$ and an even real superfunction $a=a_0+{\bf a}\in C^\infty(\R^m)\otimes \mathfrak{G}^{(ev)}_{2n}$ where $a_0$ and ${\bf a}$ are the body and nilpotent part of $a$ respectively. Then the superfunction $F(a({\bf x}))\in C^\infty(\R^m)\otimes \mathfrak{G}^{(ev)}_{2n}$ is defined, through the Taylor expansion of $F$, as
\begin{equation}\label{BosSuFun}
F(a)=F(a_0+{\bf a})=\sum_{j=0}^{n} \;\frac{{\bf a}^j}{j!}\;F^{(j)}(a_0).
\end{equation}
The easiest application of the extension (\ref{BosSuFun}) is obtained when defining arbitrary real powers of even superfunctions.    Let $a=a_0+{\bf a}\in C^\infty(\R^m)\otimes \mathfrak{G}_{2n}^{(ev)}$ be a real superfunction and $p\in \R$, then for $a_0>0$ we define 
\begin{equation}\label{PosPowBosElm}
a^p=\sum_{j=0}^{n}\;\frac{{\bf a}^j}{j!}\; (-1)^j\,(-p)_j \,a_0^{p-j}, \hspace{.5cm} \mbox{ where }  \hspace{.5cm} (q)_j=\begin{cases} 1, & j=0,\\ q(q+1)\cdots (q+j-1), & j>0, \end{cases}
\end{equation}
is the rising Pochhammer symbol. {If the numbers $q$ and $q+j$ belong to the set $\R\setminus\{0,-1,-2,\ldots\}$} we can write $\displaystyle (q)_j=\frac{\Gam(q+j)}{\Gam(q)}$. Making use of this definition of power function in superspace, we can easily prove that its basic properties still hold in this setting, i.e.\
\begin{align*}
a^pa^q&=a^{p+q}, & (ab)^p&=a^pb^p, & (a^p)^q&=a^{pq},
\end{align*}
where $p,q\in \R$ and $a,b\in C^\infty(\R^m)\otimes \mathfrak{G}_{2n}^{(ev)}$ are real superfunctions with  positive bodies.

The square of the supervector variable ${\bf x}$ yields an even superfunction with negative body, i.e.\
\[{\bf x}^2=-\sum_{j=1}^m x_j^2 + \sum_{j=1}^n x\p_{2j-1} x\p_{2j}\in C^\infty(\R^m)\otimes \mathfrak{G}_{2n}^{(ev)}.\]
Then, {\it the absolute value} of ${\bf x}$ can be defined as in the classical case by
\[|{\bf x}|=(-{\bf x}^2)^{1/2}=\left(|\underline{x}|^2- \underline{x\p}^2 \right)^{1/2}=\sum_{j=0}^n \frac{(-1)^j \underline{x\p}^{\,2j}}{j!} \, \frac{\Gam(\frac{3}{2})}{\Gam(\frac{3}{2}-j)} |\underline{x}|^{1-2j}, \;\;\mbox{ where } \;\; |\underline{x}|=\left(\sum_{j=1}^m x_j^2\right)^{\frac{1}{2}}.\]

When doubling the bosonic dimension, it is possible to define a so-called complex structure ${\bf J}$ on the algebra $\mathcal{P} \otimes\mathcal C_{2m,2n}$, see e.g.\ \cite{DS_Guz_Somm2, MR3702693}. In general, ${\bf J}$ is an algebra automorphism over $\mathcal{P} \otimes\mathcal C_{2m,2n}$ defined by
 \begin{itemize}
\item ${\bf J}$ is the identity on $\mathcal{P}$;\\[-6mm]
\item ${\bf J}(e_j)=-e_{m+j}$, ${\bf J}(e_{m+j})=e_j$, $j=1,\ldots,m$;\\
${\bf J}(e\p_{2j-1})=-e\p_{2j}$, ${\bf J}(e\p_{2j})=e\p_{2j-1}$, $j=1,\ldots,n$;\\[-6mm]
\item ${\bf J}(FG)={\bf J}(F){\bf J}(G)$ for all $F,G\in \mathcal{P} \otimes\mathcal C_{2m,2n}$.
\end{itemize}
Thence the action of ${\bf J}$ on the supervector variable ${\bf x}\in \mathcal{P} \otimes\mathcal C_{2m,2n}$ has the form
\[{\bf J}({\bf x})={\bf J}(\underline{x})+{\bf J}(\underline{x\p})=\sum_{j=1}^{m} \left(x_{m+j}e_j-x_je_{m+j}\right) +\sum_{j=1}^n \left(x\p_{2j}e\p_{2j-1}-x\p_{2j-1}e\p_{2j}\right).\]
The partial derivatives $\pa_{x_j}, \pa_{x\p_j}$ always commute with the complex structure ${\bf J}$. Then, the corresponding action of ${\bf J}$ on the super Dirac operator $\pa_{\bf x}$ can be easily seen by means of the following bosonic and fermionic twisted Dirac operators
\begin{align*}
\pa_{{\bf J}(\underline x)}&:={\bf J}(\pa_{\underline x})=\sum_{j=1}^m \left(e_j\pa_{x_{m+j}}-e_{m+j}\pa_{x_{j}}\right), & \pa_{{\bf J}(\underline x\p)}&:={\bf J}(\pa_{\underline x\p})=2\sum_{j=1}^n \left(e\p_{2j-1}\pa_{x\p_{2j-1}}+e\p_{2j}\pa_{x\p_{2j}}\right).
\end{align*}
The twisted super Dirac operator $\pa_{{\bf J}({\bf x})}$ is then defined  by
\begin{align*}
\pa_{{\bf J}({\bf x})} \cdot &:={\bf J}(\pa_{\bf x} \cdot) =\pa_{{\bf J}(\underline x\p)}\cdot -\pa_{{\bf J}(\underline x)}\cdot \,, & \cdot \pa_{{\bf J}({{\bf x}})}& :={\bf J}(\cdot \pa_{{\bf x}}) = -\cdot \pa_{{\bf J}(\underline x\p)}-\cdot \pa_{{\bf J}(\underline x)}.
\end{align*}
We recall that the actions of $\pa_{{\bf x}}$ and $\pa_{{\bf J}({\bf x})}$ on the vector variables ${\bf x}$, ${\bf J}({\bf x})$ give rise to two important defining elements: the superdimension $M=2m-2n$ and the superbivector {\bf B}, see \cite{DS_Guz_Somm2}. Indeed,
\[\pa_{\bf x}[{\bf x}]=[{\bf x}]\pa_{\bf x}=\pa_{{\bf J}({\bf x})}[{{\bf J}({\bf x})}]=[{{\bf J}({\bf x})}]\pa_{{\bf J}({\bf x})}=2m-2n,\]
while
\begin{align*}
\pa_{\underline{x}}[{\bf J}(\underline{x})] &=-[{\bf J}(\underline{x})]\pa_{\underline{x}}=-2{\bf B}_b:=-2\sum_{j=1}^m e_je_{m+j}, & \pa_{\underline{x}\p}[{\bf J}(\underline{x}\p)] &=-[{\bf J}(\underline{x}\p)]\pa_{\underline{x}\p}=-2{\bf B}_f:=-2\sum_{j=1}^{2n} e\p_{j}^{\;2},
\end{align*}
yield
\[\pa_{\bf x}[{\bf J}({\bf x})]=-[{\bf J}({\bf x})]\pa_{\bf x}=2{\bf B}=-\pa_{{\bf J}({\bf x})}[{\bf x}]=[{\bf x}]\pa_{{\bf J}({\bf x})},\]
where ${\bold B}:={\bf B}_b-{\bf B}_f=\sum_{j=1}^m e_je_{m+j}-\sum_{j=1}^{2n} e\p_j{}^2$. In addition, the following identities are easy to check
\begin{align*}
\pa^2_{{\bf x}}&= \pa^2_{{\bf J}({\bf x})}=-\Del_{2m|2n}, & \{\pa_{{\bf x}},\pa_{{\bf J}({\bf x})}\}&=0,
\end{align*}
where $\Del_{2m|2n}=\sum_{j=1}^{2m} \pa^2_{x_j}-4\sum_{j=1}^n \pa_{x\p_{2j-1}}\pa_{x\p_{2j}}$ is the super Laplace operator and the symbol $\{a,b\}:=ab+ba$ denotes the anti-commutator of a pair elements $a,b$.

The above complex structure defines the framework for so-called Hermitian Clifford analysis in the complexification $\C\mathcal{P} \otimes\mathcal C_{2m,2n}$ of $\mathcal{P} \otimes\mathcal C_{2m,2n}$, i.e.\ 
\[\C\mathcal{P} \otimes\mathcal C_{2m,2n}= \C[x_1,\ldots, x_{2m}]\otimes \mathfrak{G}_{2n} \otimes \mathcal C_{2m,2n}=\mathcal{P} \otimes\mathcal C_{2m,2n} \oplus i \; \mathcal{P} \otimes\mathcal C_{2m,2n},\]
where $\C[x_1,\ldots, x_{2m}]$ denotes the complexification of $\R[x_1,\ldots, x_{2m}]$. The Hermitian and complex conjugations are defined on $\C\mathcal{P} \otimes\mathcal C_{2m,2n}$ by 
\begin{align*}
(a+ib)^\dagger&=\ba{a}-i\ba{b}, & (a+ib)^c&={a}-i{b}, & a,b&\in \mathcal{P} \otimes\mathcal C_{2m,2n},
\end{align*}
respectively. Here the bar notation stands for the Clifford conjugation in $ \mathcal{P} \otimes\mathcal C_{2m,2n}$, i.e. the linear mapping satisfying 
\[a \;e_{j_1}\ldots e_{j_k} e\p_{\ell_1}\ldots e\p_{\ell_s}=a \;\ba{e_{j_1}\ldots e_{j_k} e\p_{\ell_1}\ldots e\p_{\ell_s}}, \hspace{.3cm} a\in\mathcal P;  \]   
where
\[\ba{e_{j_1}\ldots e_{j_k} e\p_{\ell_1}\ldots e\p_{\ell_s}}=(-1)^{k+\frac{s(s+1)}{2}}e\p_{\ell_s}\ldots e\p_{\ell_1} e_{j_k}\ldots e_{j_1}.  \]

As in the classical case, the two projection operators $\frac{1}{2}\left(1\pm i{\bf J}\right)$ will produce the main objects of the Hermitian setting by acting upon the corresponding objects in $\mathcal{P} \otimes\mathcal C_{2m,2n}$. In first place, the so-called Witt basis elements in superspace are obtained through the action of $\pm\frac{1}{2}\left(1\pm i{\bf J}\right)$ on the orthogonal and symplectic generators $e_j$, $e\p_j$, i.e.\
\begin{align*}
\mathfrak f_j &=\frac{1}{2} (e_j-ie_{m+j}), & \mathfrak f_j^\dag&=-\frac{1}{2} (e_j+ie_{m+j}), & \mathfrak f_j\p&=\frac{1}{2} (e\p_{2j-1}-ie\p_{2j}), &\mathfrak f_j\p^\dag&=-\frac{1}{2} (e\p_{2j-1}+ie\p_{2j}).
\end{align*}
These Witt basis elements submit to the multiplication rules
\[
\begin{cases}
\mathfrak f_j \mathfrak f_k+\mathfrak f_k \mathfrak f_j=0,\\
\mathfrak f_j^\dag \mathfrak f_k^\dag+\mathfrak f_k^\dag \mathfrak f_j^\dag=0,\\
\mathfrak f_j \mathfrak f_k^\dag+\mathfrak f_k^\dag \mathfrak f_j=\del_{j,k},
\end{cases}
\hspace{.3cm}
\begin{cases}
\mathfrak f_j\p \,\mathfrak f_k\p - \mathfrak f_k\p\, \mathfrak f_j\p=0,\\
\mathfrak f_j\p^\dag \,\mathfrak f_k\p^\dag-\mathfrak f_k\p^\dag \,\mathfrak f_j\p^\dag=0,\\
\mathfrak f_j\p\, \mathfrak f_k\p^\dag-\mathfrak f_k\p^\dag \,\mathfrak f_j\p=-\frac{i}{2}\del_{j,k},
\end{cases}
\hspace{.3cm}
\begin{cases}
\mathfrak f_j\mathfrak f_k\p + \mathfrak f_k\p\, \mathfrak f_j=0,\\
\mathfrak f_j \mathfrak f_k\p^\dag+\mathfrak f_k\p^\dag \,\mathfrak f_j=0,
\end{cases}
\hspace{.3cm}
\begin{cases}
\mathfrak f_j^\dag \mathfrak f_k\p + \mathfrak f_k\p\, \mathfrak f_j^\dag =0,\\
\mathfrak f_j^\dag  \mathfrak f_k\p^\dag+\mathfrak f_k\p^\dag \,\mathfrak f_j^\dag =0.
\end{cases}
\]
The actions of the projection operators on the supervector variable ${\bf x}$ produce the Hermitian supervector variable ${\bf Z}$ and its Hermitian conjugate ${\bf Z}^\dag$
\begin{align*}
{\bf Z}&=\phantom{-}\frac{1}{2}\Big({\bf x}+i{\bf J}({\bf x})\Big)=\underline{Z}+\underline{Z^{\p}}\phantom{{}^\dagger}=\sum_{j=1}^m z_j{\mathfrak f}_j+\sum_{j=1}^n z\p_j{{\mathfrak f}_j}\p,\\
{\bf Z}^\dagger &= -\frac{1}{2}\Big({\bf x}-i{\bf J}({\bf x})\Big)=\underline{Z}^\dagger+\underline{{Z^{\p}}}^\dagger= \sum_{j=1}^m z_j^c{\mathfrak f}_j^\dagger+\sum_{j=1}^n z\p_j^c\, {{\mathfrak f_j}\p}^\dag ,
\end{align*}
where the bosonic Hermitian vector variables $\underline{Z}$,$\underline{Z}^\dagger$ and the fermionic ones $\underline{Z}^{\p}$,${\underline{Z}^{\p}}^\dagger$ are given by
\begin{align*}
\underline{Z}&=\frac{1}{2}\Big(\underline{x}+i{\bf J}(\underline{x})\Big)=\sum_{j=1}^m z_j{\mathfrak f}_j, &
\underline{Z}^\dagger&=-\frac{1}{2}\Big(\underline{x}-i{\bf J}(\underline{x})\Big)=\sum_{j=1}^m z_j^c{\mathfrak f}_j^\dagger,\\
\underline{Z}^{\p}&= \frac{1}{2}\Big(\underline{x\p}+i{\bf J}(\underline{x\p})\Big)=\sum_{j=1}^n z\p_j{{\mathfrak f}_j}\p, & {\underline{Z}^{\p}}^\dagger&= -\frac{1}{2}\Big(\underline{x\p}-i{\bf J}(\underline{x\p})\Big)=\sum_{j=1}^n z\p_j^c\, {{\mathfrak f_j}\p}^\dag. 
\end{align*}
Here we have introduced the commuting and anti-commuting complex variables $z_j=x_j+ix_{m+j}$, $z\p_j=x\p_{2j-1}+ix\p_{2j}$ together with their complex conjugates $z_j^c=x_j-ix_{m+j}$, $z\p_j^c=x\p_{2j-1}-ix\p_{2j}$.

Finally, the Hermitian Dirac operators $\pa_{{\bf Z}}$, $\pa_{{\bf Z}^\dag}$ are derived from $\pa_{\bf x}$ by means of the relations
\begin{align}\label{FormPaXtoPaZ}
\pa_{{\bf Z}} &:= \frac{1}{4}\left(\pa_{\bf x}-i \pa_{{\bf J}({\bf x})}\right), & \pa_{{\bf Z}^\dag} &:= -\frac{1}{4}\left(\pa_{\bf x}+i \pa_{{\bf J}({\bf x})}\right),
\end{align}
which are valid for both left and right actions of the operators $\pa_{{\bf Z}}$, $\pa_{{\bf Z}^\dag}$. These actions can be re-written as 
\begin{align*}
\begin{cases}
\pa_{\bf Z} \cdot =\pa_{{\underline{Z}^{\p}}}\cdot+\pa_{\underline{Z}}\cdot, \\ 
\cdot\pa_{\bf Z}  =-\cdot\pa_{{\underline{Z}^{\p}}}+\cdot\pa_{\underline{Z}},
\end{cases}
&&
\begin{cases}
\pa_{{\bf Z}^\dag} \cdot = \pa_{{\underline{Z}^{\p}}^\dagger}\cdot+\pa_{\underline{Z}^\dag}\cdot,\\
\cdot \pa_{{\bf Z}^\dag}= - \cdot \pa_{{\underline{Z}^{\p}}^\dagger}+\cdot \pa_{\underline{Z}^\dag}, 
\end{cases}
\end{align*}
where we have introduced the bosonic and fermionic Hermitian Dirac operators 
\begin{align*}
\pa_{\underline{Z}} &=-\frac{1}{4} \big(\pa_{\underline{x}}-i\pa_{{\bf J}(\underline{x})}\big)=\sum_{j=1}^m \f_j^\dag \,\pa_{z_j}, &  \pa_{\underline{Z}^\dag}&=  \frac{1}{4} \big(\pa_{\underline{x}}+i\pa_{{\bf J}(\underline{x})}\big)=\sum_{j=1}^m \f_j \,\pa_{z_j^c},  \\
\pa_{{\underline{Z}^{\p}}}& = \frac{1}{4} \big(\pa_{\underline{x\p}} -i\pa_{{\bf J}(\underline{x\p})}\big)=2i\sum_{j=1}^n \f_j\p^\dag \, \pa_{z\p_j}, &
 \pa_{{\underline{Z}^{\p}}^\dag}&= -\frac{1}{4} \big(\pa_{\underline{x\p}} +i\pa_{{\bf J}(\underline{x\p})}\big)=-2i\sum_{j=1}^n \f_j\p \, \pa_{z\p_j^c};
\end{align*}
involving the classical Cauchy-Riemann operators $\pa_{z_j}=\frac{1}{2}(\pa_{x_j}-i\pa_{x_{m+j}})$, $\pa_{z\p_j}=\frac{1}{2}(\pa_{x\p_{2j-1}}-i\pa_{x\p_{2j}})$  and their conjugates $\pa_{z_j^c}=\frac{1}{2}(\pa_{x_j}+i\pa_{x_{m+j}})$, 
$\pa_{z\p_j^c}=\frac{1}{2}(\pa_{x\p_{2j-1}}+i\pa_{x\p_{2j}})$ with respect to the variables $z_j$ and $z\p_j$.

As was the case with $\pa_{\bf x}$, the notion of super monogenicity may be naturally associated to $\pa_{{\bf J}({\bf x})}$ as well. Then a function $F\in C^1(\Om)\otimes \mathfrak{G}_{2n} \otimes \mathcal{C}_{2m,2n}$, where $\Om\inc \R^{2m}$ is an open set, is called a (left) super Hermitian monogenic (or sh-monogenic) function if it satisfies the system 
\[\pa_{\bf x}[F]=0=\pa_{{\bf J}({\bf x})}[F]\]
or equivalently, the system
\[\pa_{{\bf Z}}[F]=0=\pa_{{\bf Z}^\dag}[F].\]
It can be easily proven that 
\begin{align*}
\left({\bf Z}\right)^2&=\left({\bf Z}^\dagger \right)^2=0,  & \left(\pa_{{\bf Z}} \right)^2&= \left(\pa_{{\bf Z}^\dag} \right)^2=0, & \Del_{2m|2n}&=4\left\{\pa_{{\bf Z}}, \pa_{{\bf Z}^\dag}\right\}.
\end{align*}
Moreover, if we define $|{\bf Z}|^2=|{\bf Z}^\dag|^2:=\left\{{{\bf Z}}, {{\bf Z}^\dag}\right\}$ one immediately has
\[|{\bf Z}|^2=|{\bf Z}^\dag|^2=\sum_{j=1}^m z_jz_j^c-\frac{i}{2}\sum_{j=1}^n z\p_jz\p_j^c=\sum_{j=1}^{2m} x_j^2 - \sum_{j=1}^n x\p_{2j-1} x\p_{2j}=|{\bf x}|^2=|{\bf J}({\bf x})|^2.\]

\section{Distributions and integration in superspace}
In this section we recall some important aspects of distributions in superanalysis, and their use in integration over general surfaces and domains in superspace. For the convenience of the reader we repeat the relevant material from \cite{Guz_Somm5}, thus making our exposition self-contained.

\subsection{Distributions in superanalysis}
Let $\mathcal{D}^\prime$ be the space of Schwartz distributions,  i.e. the space of generalized functions on the space $C^\infty_0(\R^m)$ of complex valued $C^\infty$-functions with compact support. As usual, the notation
\begin{equation}\label{DistEv}
\int_{\R^m} \al f \; dV_{\underline x}= \langle\al,f\rangle,
\end{equation}
where $dV_{\underline{x}}=d{x_1}\cdots d{x_m}$ is the classical $m$-volume element, is used for the evaluation of the distribution $\al\in \mathcal{D}^\prime$ on the test function $f\in C^\infty_0(\R^m)$.

Let  $\mathcal{E}^\prime$ be the space of  generalized functions on the space  $C^\infty (\R^m)$ of $C^\infty$-functions in $\R^m$ (with arbitrary support). We recall that $\mathcal{E}^\prime$ is exactly the subspace of all compactly supported distributions in $\mathcal{D}^\prime$. Indeed, every distribution in $\mathcal{E}^\prime\inc \mathcal{D}^\prime$ has compact support and vice-versa, every distribution in $\mathcal{D}^\prime$ with compact support can be uniquely extended to a distribution in $\mathcal{E}^\prime$, see \cite{MR0208364} for more details. This means that, for every $\al\in \mathcal{E}^\prime$, evaluations of the form (\ref{DistEv}) extend to $C^\infty (\R^m)$ (instead of $C^\infty_0 (\R^m)$).

The space of superdistributions  $\mathcal{D}^\prime \otimes \mathfrak{G}_{2n}$ {then} is defined by all the elements of the form
\begin{equation}\label{SupDistForm}
\al=\sum_{A\inc\{1,\ldots, 2n\}} \al_A \underline{x}\p_A, \hspace{.5cm} \mbox{ with } \al_A\in \mathcal{D}^\prime.
\end{equation}
Similarly, the subspace $\mathcal{E}^\prime \otimes \mathfrak{G}_{2n}$ is composed by all the elements of the form (\ref{SupDistForm}) but with $\al_A\in \mathcal{E}^\prime$.

The {analogue} of the integral $\int_{\R^m}\; dV_{\underline x}$ in superspace  is given by
\[\int_{R^{m|2n}}=\int_{\R^m} dV_{\underline{x}} \int_B=\int_B \int_{\R^m} dV_{\underline{x}},\]
where the bosonic integration is the usual real integration and the integral over fermionic variables is given  by the so-called Berezin integral (see \cite{Berezin:1987:ISA:38130}), defined by
\[\int_B=\pi^{-n} \, \pa_{x\p_{2n}}\cdots \pa_{x\p_{1}}=\frac{(-1)^n \pi^{-n}}{4^n n!} \pa_{\underline{x}\p}^{2n}.\]
This {enables us} to define the action of a superdistribution $\al\in \mathcal{D}^\prime \otimes \mathfrak{G}_{2n}$ (resp. $\al\in \mathcal{E}^\prime \otimes \mathfrak{G}_{2n}$) on a test superfunction $F\in C^\infty_0(\R^m) \otimes \mathfrak{G}_{2n}$ (resp. $F\in C^\infty(\R^m) \otimes \mathfrak{G}_{2n}$) by
\[\int_{\R^{m|2n}} \al F := \sum_{A,B\inc\{1,\ldots, 2n\}} \langle \al_A, f_B \rangle \int_B   \underline{x}\p_A\, \underline{x}\p_B. \]
As in the classical case, we say that the superdistribution $\al\in \mathcal{D}^\prime \otimes \mathfrak{G}_{2n}$ {\it vanishes in the open set} $\Om\inc \R^m$ if $\int_{\R^{m|2n}} \al F=0$ for every $F\in C^\infty_0(\R^m) \otimes \mathfrak{G}_{2n}$ whose real support is contained in $\Om$. Same way, the {\it support $supp\; \al$ of} $\al\in \mathcal{D}^\prime \otimes \mathfrak{G}_{2n}$ is defined as the complement of the largest open subset of $\R^m$ on which $\al$ vanishes. Hence, it can be easily seen that $supp\; \al=\bigcup_{A\inc\{1,\ldots, 2n\}} supp\; \al_A$. This means that  $\mathcal{E}^\prime \otimes \mathfrak{G}_{2n}$ is the subspace of all compactly supported superdistributions in $\mathcal{D}^\prime \otimes \mathfrak{G}_{2n}$.

We now define the multiplication of distributions with disjoint singular supports. We first recall that the singular support $sing \; supp \, \al$ of the distribution $\al\in\mathcal{D}'$ is defined by {the statement} that $\underline{x}\notin sing \; supp \, \al$ if and only if there exists a neighborhood $U_{\underline{x}}$ of $\underline{x}\in\R^m$ such that the restriction of $\al$ to $U_{\underline{x}}$ is a smooth function. It is readily seen that $sing \; supp \, \al \inc supp \, \al$.
\begin{defi}[Multiplication of distributions, \cite{MR3270560, MR1996773}]\label{MultDistDisSingSupp}
Consider two distributions $\al,\be\in \mathcal{D}'$ such that $sing \; supp \, \al \cap sing \; supp \, \be=\emptyset$. The product of distributions $\al\be$ is well defined {by} the formula
\begin{equation}\label{DistProd}
\langle \al\be, \phi\rangle= \langle \al, \be \chi \phi\rangle + \langle \be, \al (1-\chi) \phi\rangle, \hspace{1cm} \phi\in C^{\infty}(\R^m),
\end{equation}
where $\chi\in C^\infty(\R^m)$  is equal to zero in a neighborhood of $sing \; supp \, \be$ and equal to one in a neighborhood of $sing \; supp \, \al$.
\end{defi}
\begin{remark}
In general, the product of two distributions in $\mathcal{D}^\prime$ is defined if their wave front sets satisfy the so-called {\it Hörmanders condition}. See \cite[p.~267]{MR1996773} and \cite{MR3270560} for more details.
\end{remark}
It is easily seen that if $\al,\be\in \mathcal{D}'$ vanish in $\Om\inc\R^m$ then the product $\al\be$ vanishes in $\Om$ as well. Hence $supp\, \al\be \inc supp\,\al \cup supp\, \be$.
As a consequence, if $\al$ and $\be$ have compact supports (i.e.\ $\al,\be\in  \mathcal{E}'$) then $\al\be$ also has compact support (i.e.\ $\al\be\in  \mathcal{E}'$). The product (\ref{DistProd}) is associative, commutative and satisfies the Leibniz rule, see \cite{MR3270560, MR1996773}. 

The notion of singular support can be extended to distributions $\al\in \mathcal{D}'\otimes \mathfrak{G}_{2n}$ by {the statement} that $\underline{x}\notin sing \; supp \, \al$ if and only if there exists a neighborhood $U_{\underline{x}}$ of $\underline{x}\in\R^m$ such that the restriction of $\al$  to $U_{\underline{x}}$ belongs to $C^\infty(U_{\underline x}) \otimes \mathfrak{G}_{2n}$.  {In} this way we obtain for every $\al\in \mathcal{D}'\otimes \mathfrak{G}_{2n}$ of the form (\ref{SupDistForm}) that $sing \; supp \, \al =\bigcup_{A\inc\{1,\ldots, n\}} sing \; supp \, \al_A$.
In the same way, we define the product of superdistributions $\al, \be\in \mathcal{D}'\otimes \mathfrak{G}_{2n}$ with $sing \; supp \, \al \cap sing \; supp \, \be=\emptyset$ by
\begin{equation}\label{prodSupDist}
\al\be = \ \sum_{A,B\inc\{1,\ldots, n\}} \al_A\be_B\; \underline{x}\p_A\, \underline{x}\p_B,
\end{equation}
where the distribution $\al_A\be_B$ is {to be} understood in the sense of (\ref{DistProd}).

\subsection{Domain and surface integrals in superspace}
The notions of general domain and surface integrals in superspace were introduced in \cite{Guz_Somm5} by means of the corresponding extensions of the Dirac and Heaviside distributions. In  \cite{MR2539324}, these distributions were introduced for some particular cases corresponding to the supersphere.

From now on, we use the notation $p\N+q:=\{pk+q:k\in\N\}$ where $p,q\in \Z$ and $\N:=\{1,2,\ldots\}$ is the set of natural numbers.

Consider an even real superfunction $g=g_0+{\bf g}\in C^\infty(\R^m)\otimes \mathfrak{G}_{2n}^{(ev)}$ such that $\pa_{\underline x}[g_0]\neq0$ on the surface $g_0^{-1}(0)$. The distribution $\del^{(k)}(g)$ is defined as the Taylor series  
\[\del^{(k)}(g) = \sum_{j=0}^{n} \frac{{\bf g}^j}{j!} \, \del^{(k+j)}(g_0), \hspace{1cm} k\in\N-2:=\{-1, 0, 1, 2\ldots\}.
\]
The particular case $k=-1$ provides the expression for the antiderivative of $\del$, i.e. the Heaviside distribution $H=\del^{(-1)}$ given by
\begin{equation}\label{HaevSupS}
H(g)=H(g_0)+\sum_{j=1}^{n} \frac{{\bf g}^j}{j!} \, \del^{(j-1)}(g_0), \hspace{1cm} \mbox{ where } \hspace{.5cm}H(g_0)=\begin{cases}  1, & g_0\geq 0,\\ 0, & g_0<0.\end{cases}
\end{equation}

%
In this way, one can consider domains in superspace $\Om_{m|2n}$ given {by} characteristic functions of the form $H(-g)$. The superfunction $g$ is called the defining {\it phase function} of the domain $\Om_{m|2n}$. Observe that $\Om_{m|2n}$ plays the same r\^ole in superanalysis as its associated real domain $\Om_{m|0}:=\{\underline{x}\in\R^m: g_0(\underline{x})<0\}$ in classical analysis.
\begin{defi}[Domain integration in superspace, \cite{Guz_Somm5}]\label{DomIntSupS}
Let $\Om_{m|2n}$ be a domain in superspace (defined as before) satisfying the following two conditions:
\begin{itemize}
\item the associated real domain $\Om_{m|0}$ has compact closure;
\item the body $g_0$ of the defining  phase function $g$ is such that $\pa_{\underline x}[g_0]\neq0$ on $g_0^{-1}(0)$.
\end{itemize}
The integral over $\Om_{m|2n}$ {then} is defined as the functional $\int_{\Om_{m|2n}}:C^{n-1}\left(\ba{\Om_{m|0}}\right)\otimes \mathfrak{G}_{2n}\fd \C$ given by
\begin{equation}\label{DomSupInt}
\int_{\Om_{m|2n}} F= \int_{\R^{m|2n}} H(-g) F, \hspace{2cm} F\in C^{n-1}\left(\ba{\Om_{m|0}}\right)\otimes \mathfrak{G}_{2n}.
\end{equation}
\end{defi}
The evaluation of the expression on (\ref{DomSupInt}) requires the integration of smooth functions on the real domain $\Om_{m|0}$ {to be possible}. This is guaranteed by the first condition imposed on the super-domain $\Om_{m|2n}$. On the other hand, if the nilpotent part ${\bf g}$ of $g$ is not identically zero, the above definition {also} involves the action of the Dirac distribution {on} $g_0$, see (\ref{HaevSupS}). For that reason, we restrict our analysis to the case {where} $\pa_{\underline x}[g_0]\neq0$ on $g_0^{-1}(0)$, {in order to ensure that this action } is well defined.

The most simple examples {for illustrating} the use of Definition \ref{DomIntSupS} {correspond to the cases} $g=g_0$ or $g=-{\bf x}^2-R^2$; i.e. integration over real domains or over the superball respectively. The integral (\ref{DomSupInt}) {then} is given by
\[\int_\Om \int_B F dV_{\underline x},\hspace{.5cm} \mbox{ {and} } \hspace{.5cm}\int_{\R^{m|2n}} H({\bf x}^2+R^2) F,\]
respectively. These are the two particular cases that have been treated in the literature, see \cite{MR2539324, MR2521367}. 

In real analysis {the choice of} the real-valued function $g_0$ {defining} a certain domain $\Om_{m|0}\inc\R^m$ {is not unique}. Indeed, for every function $h_0\in C(\R^m)$ {with} $h_0>0$, the function $h_0g_0$ defines the same domain as $g_0$, i.e.\ $\Om_{m|0}=\{\underline{x}\in \R^m: g_0(\underline{x})<0\}=\{\underline{x}\in \R^m: h_0(\underline{x})g_0(\underline{x})<0\}$.
However, integration over $\Om_{m|0}$ remains  independent of the choice of the function $g_0$ that defines $\Om_{m|0}$. It suffices to note that $H(-g_0)=H(-h_0g_0)$ for $h_0>0$. This property remains valid in superspace. Indeed, in  \cite{Guz_Somm5} it was proven that for  any other phase function $hg$, where $h=h_0+{\bf h}\in C^\infty(\R^m)\otimes \mathfrak{G}_{2n}^{(ev)}$ has a positive body $h_0>0$, one has $H(hg)=H(g)$. 

{Similarly to} the case of super domains, we define a surface $\Gam_{m-1|2n}$ in superspace by means of $\del(g)$ where $g({\bf x})=g_0(\underline{x})+{\bf g}(\underline{x}, \underline{x}\p)\in C^\infty(\R^m)\otimes \mathfrak{G}_{2n}^{(ev)}$ is a phase function. If $\Om_{m|2n}$ is the super domain associated to $g$ as in Definition \ref{DomIntSupS}, then we say that $\Gam_{m-1|2n}$ is the boundary of $\Om_{m|2n}$ and denote it by $\Gam_{m-1|2n}:=\pa\Om_{m|2n}$. This way, $\Gam_{m-1|2n}$ plays the same r\^ole in superspace as its real surface $\Gam_{m-1|0}:=\pa\Om_{m|0}=\{\underline{x}\in \R^m: g_0(\underline{x})=0\}$ in classical analysis.
\begin{defi}\label{SurfIntSupS}
Let  $\Gam_{m-1|2n}$ be a surface in superspace (defined as before) satisfying the following two conditions:
\begin{itemize}
\item the associated real surface $\Gam_{m-1|0}\inc \R^m$ is a compact set;
\item the body $g_0$ of the defining  phase function $g$ is such that $\pa_{\underline x}[g_0]\neq0$ on $g_0^{-1}(0)$.
\end{itemize}
The non-oriented and oriented surface integrals over $\Gam_{m-1|2n}$ {then} are defined as the following functionals  on $C^{n}\left(\Gam_{m-1|0}\right)\otimes \mathfrak{G}_{2n}$\begin{equation}\label{SurfIntNO}
\int_{\Gam_{m-1|2n}} F = \int_{\R^{m|2n}} \del(g) \, \big|\pa_{{\bf x}}[g]\big| \, F,\hspace{.8cm} \int_{\Gam_{m-1|2n}} \sigma_{\bf x}\, F = \int_{\R^{m|2n}} \del(g) \, \pa_{{\bf x}}[g] \, F,
\end{equation}
respectively.
\end{defi}
\begin{remark}
Since $\del(g)=\del(-g)$,  the sign of the superfunction $g$ does not {play a r\^ole} in the non-oriented case, see \cite{Guz_Somm5}.
\end{remark}
When $g=g_0$, the integrals (\ref{SurfIntNO}) reduce to the product of the classical real surface integration and the Berezin integral, i.e.,
\begin{equation}\label{RealTimBer}
\int_{\Gam_{m-1|0}} \int_B F(\underline{w}, \underline{x}\p)\, dS_{\underline{w}}, \hspace{.8cm} \int_{\Gam_{m-1|0}} \int_B n(\underline{w})F(\underline{w}, \underline{x}\p)\, dS_{\underline{w}},
\end{equation}
while for $g=-{\bf x}^2-R^2$ one obtains the integral over the supersphere of radius $R>0$ given in \cite{MR2539324}, i.e.\
\[2 R \int_{\R^{m|2n}}  \del({\bf x}^2+R^2)\, F({\bf x}) .\]

As {for} domain integrals,  {in \cite{Guz_Somm5}} it was proven that Definition \ref{SurfIntSupS} does not depend on the choice of the defining phase function $g$ for the surface $\Gam_{m-1|2n}$. Indeed, for  any other phase function $hg$, where $h=h_0+{\bf h}\in C^\infty(\R^m)\otimes \mathfrak{G}_{2n}^{(ev)}$ has a positive body $h_0>0$, one has
\begin{align*}
\del(hg)\pa_{\bf x}[hg]&=\del(g)\pa_{\bf x}[g], & \del(hg) |\pa_{\bf x}[hg]|&=\del(g)|\pa_{\bf x}[g]|.
\end{align*}

When considering $2m$ bosonic variables and an open region $\Om\inc \R^{2m}$, the theorem of Stokes may be formulated as follows. 
\begin{teo}[{\bf Distributional Stokes Theorem, \cite{MR2539324}}]{}\label{StokTe}
Let $F,G\in C^\infty(\Om)\otimes \mathfrak{G}_{2n} \otimes \mathcal C_{m,2n}$ and $\al\in \mathcal{E}' \otimes \mathfrak{G}_{2n}^{(ev)}$ a distribution with compact support such that $supp\, \al\inc \Om$. Then,
\begin{align}
\int_{\R^{2m|2n}} \left(F\pa_{\bf x}\right)\al G+F\al \left(\pa_{\bf x}G\right)&=-\int_{\R^{2m|2n}} F\left(\pa_{\bf x}\al\right)G, \label{Stok1} \\
\int_{\R^{2m|2n}} \left(F\pa_{{\bf J}({\bf x})}\right)\al G+F\al \left(\pa_{{\bf J}({\bf x})}G\right)&=-\int_{\R^{2m|2n}} F\left(\pa_{{\bf J}({\bf x})}\al\right)G. \label{Stok2}
\end{align}
\end{teo}

Formula (\ref{Stok1}) was obtained in \cite{MR2539324}, and formula (\ref{Stok2}) is easily derived by letting the complex structure ${\bf J}$ act on both sides of (\ref{Stok1}).

In particular, for $\al=H(-g)$  
it holds that (see \cite{Guz_Somm5})
\begin{align*}
\pa_{\bf x}[H(-g)]&=-\del(g) \pa_{\bf x}[g], & \pa_{{\bf J}({\bf x})}[H(-g)]&=-\del(g) \pa_{{\bf J}({\bf x})}[g].
\end{align*}
Hence, if we assume the set $\{g_0\leq 0\}:=\{\underline{w}\in \R^{2m}: g_0(\underline{w})\leq 0\}$ to be compact, we obtain a Stokes formula in superspace compatible with the notions of domain and surface integrals that we have introduced above.
\begin{cor}\label{StokCor}
Let $g=g_0+{\bf g}\in C^\infty(\R^{2m})\otimes \mathfrak{G}_{2n}^{(ev)}$ be a phase function such that $\{g_0\leq 0\}$ is compact and $\pa_{\underline x}[g_0]\neq0$ on $g_0^{-1}(0)$. Then, for $F,G\in C^\infty(\Om)\otimes \mathfrak{G}_{2n} \otimes \mathcal C_{m,2n}$ such that $\{g_0\leq 0\}\inc \Om$ one has
\begin{align*}
\int_{\R^{2m|2n}} H(-g )\left[\left(F\pa_{\bf x}\right) G+F \left(\pa_{\bf x}G\right)\right]& = \int_{\R^{2m|2n}} F\del(g) \pa_{\bf x}[g] G,\\ 
\int_{\R^{2m|2n}} H(-g )\left[\left(F\pa_{{\bf J}({\bf x})}\right) G+F \left(\pa_{{\bf J}({\bf x})}G\right)\right]& = \int_{\R^{2m|2n}} F\del(g) \pa_{{\bf J}({\bf x})}[g] G. 
\end{align*}
\end{cor} 
Using (\ref{FormPaXtoPaZ}) we easily obtain the following results which are equivalent to Theorem \ref{StokTe} and Corollary \ref{StokCor}.
\begin{cor}\label{Her_stokThemCor2}
The following formulae hold under the same conditions {as in} Theorem \ref{StokTe}
\begin{align}
\int_{\R^{2m|2n}} \left(F\pa_{\bf Z}\right)\al G+F\al \left(\pa_{\bf Z}G\right)&=-\int_{\R^{2m|2n}} F\left(\pa_{\bf Z}\al\right)G, \label{Stok11} \\
\int_{\R^{2m|2n}} \left(F\pa_{{\bf Z}^\dag}\right)\al G+F\al \left(\pa_{{\bf Z}^\dag}G\right)&=-\int_{\R^{2m|2n}} F\left(\pa_{{\bf Z}^\dag}\al\right)G. \label{Stok22}
\end{align}
\end{cor}
\begin{cor}
The following formulae hold under the same conditions {as in} Corollary \ref{StokCor}
\begin{align*}
\int_{\R^{2m|2n}} H(-g )\left[\left(F\pa_{\bf Z}\right) G+F \left(\pa_{\bf Z}G\right)\right]& = \int_{\R^{2m|2n}} F\del(g) \pa_{\bf Z}[g] G, \\
\int_{\R^{2m|2n}} H(-g )\left[\left(F\pa_{{\bf Z}^\dag}\right) G+F \left(\pa_{{\bf Z}^\dag}G\right)\right]& = \int_{\R^{2m|2n}} F\del(g) \pa_{{\bf Z}^\dag}[g] G. 
\end{align*}
\end{cor}

\section{Fundamental solutions for $\pa_{\bf x}$ and $\pa_{{\bf J}({\bf x})}$}
In \cite{MR2386499}, the fundamental solution of the super Dirac operator $\pa_{\bf x}$ for the general superdimension $M=m-2n$ was calculated to be,
\begin{equation}\label{FundSolSupDirOp}
\nu_1^{m|2n}=\pi^n \sum_{k=0}^{n-1} \frac{2^{2k+1} k!}{(n-k-1)!}\, \fhi_{2k+2}^{m|0} \underline{x}\p^{\, 2n-2k-1} - \pi^n \sum_{k=0}^n \frac{2^{2k} k!}{(n-k)!}\, \fhi_{2k+1}^{m|0} \underline{x}\p^{\, 2n-2k}, 
\end{equation}
where $ \fhi_j^{m|0}$ is the fundamental solution of $\pa_{\underline x}^j$. Observe that $\nu_1^{m|0}=-\fhi_1^{m|0}$. The superdistribution $\nu_1^{m|2n}$ satisfies 
\[\pa_{\bf x}\nu_1^{m|2n}({\bf x})=\del(\underline{x}) \frac{\pi^n}{n!} \underline{x}\p^{\, 2n} = \del({\bf x})= \nu_1^{m|2n}({\bf x}) \pa_{\bf x},\]
where $\del({\bf x})=\del(\underline{x}) \frac{\pi^n}{n!} \underline{x}\p^{\, 2n} $ defines the Dirac distribution on the supervector variable ${\bf x}$ and $\del(\underline{x})=\del(x_1)\cdots \del(x_{m})$ is the $m$-dimensional real Dirac distribution. It is easily seen that
\begin{equation}\label{DelProp}
\langle \del({\bf x}-{\bf y}), G({\bf x})\rangle= \int_{\R^{m|2n}} \del({\bf x}-{\bf y})G({\bf x})=G({\bf y}) \mbox{ or equivalently, } \displaystyle\begin{cases} \displaystyle \int_{\R^{m}}  \del(\underline{x}- \underline{y}) G_A(\underline{x}) \, dV_{\underline{x}}=G_A(\underline{y}), \\ \\ \displaystyle\frac{\pi^n}{n!} \int_B (\underline{x}\p-\underline{y}\p)^{2n}\, \underline{x}\p_A= \underline{y}\p_A,\end{cases}
\end{equation}
where ${\bf y}=\underline{y}+\underline{y}\p$ and $G\in C^\infty(U_{\underline y})\otimes \mathfrak{G}_{2n}$ with $U_{\underline y}\inc \R^{m}$ being a neighborhood of ${\underline y}$.
\begin{teo}
If  $M=m-2n \notin -2\N+2$, the fundamental solution  $\nu_1^{m|2n}$ has the form
\[\nu_1^{m|2n}({\bf x})=\frac{1}{|\Sa^{m-1|2n}|} \frac{{\bf x}}{|{\bf x}|^M},\]
where $|\Sa^{m-1|2n}|=\frac{2\pi^{M/2}}{\Gamma\left(\frac{M}{2}\right)}$ is the surface area of the unit supersphere, see \cite{Guz_Somm5}.
\end{teo}
\pf
We first recall that the fundamental solution of $\Del_{m}^k$, where $\Del_m=\sum_{j=1}^m \pa_{x_j}^2$ is the Laplace operator in $m$ bosonic dimensions, is given by 
\begin{equation}\label{FunPolyHarSol}
\nu_{2k}^{m|0}(\underline{x})=\frac{|\underline{x}|^{2k-n}}{|\Sa^{m-1}|\, 2^{k-1} \,(k-1)! \,\prod_{\ell=1}^k (2\ell-m)}=\frac{(-1)^k \Gam\left(\frac{m}{2}-k\right)}{2^{2k} \, \pi^{m/2} \, \Gam(k)} \, \frac{|\underline{x}|^{2k}}{|\underline{x}|^m},
\end{equation}
if $m-2k\notin -2\N+2$, see \cite{MR745128}. Then the above formula can be used for every $k\leq n$, since the condition $m-2n\notin -2\N+2$ directly implies $m-2k\notin -2\N+2$. Indeed, it suffices to observe that $m-2k=m-2n+2(n-k)$ {with} $n-k\geq 0$.

\noindent Since $\Del_{m}^{k+1}=(-1)^{k+1} \pa_{\underline{x}}^{2k+2}$ we can write
\begin{align}
\fhi_{2k+2}^{m|0} &= (-1)^{k+1} \nu_{2k+2}^{m|0} = \frac{\Gam\left(\frac{m}{2}-k-1\right)}{2^{2k+2} \, \pi^{m/2} \, \Gam(k+1)} \, \frac{|\underline{x}|^{2k+2}}{|\underline{x}|^m}, & k&=0, 1, \dots, n-1, \label{PolyMonFundSolEven} \\
\fhi_{2k+1}^{m|0} &= \pa_{\underline{x}}\left[\fhi_{2k+2}^{m|0}\right] = - \frac{\Gam\left(\frac{m}{2}-k\right)}{2^{2k+1} \, \pi^{m/2} \, \Gam(k+1)} \, \frac{\underline{x} \, |\underline{x}|^{2k}}{|\underline{x}|^m}, & k&=0, 1, \dots, n-1. \label{PolyMonFundSolOdd}
\end{align}
It is easily seen that (\ref{PolyMonFundSolOdd}) also holds for $k=n$. Indeed, {writing}
\[\fhi_{2n+1}^{m|0} = - \frac{\Gam\left(\frac{m}{2}-n\right)}{2^{2n+1} \, \pi^{m/2} \, \Gam(n+1)} \, \frac{\underline{x} \, |\underline{x}|^{2n}}{|\underline{x}|^m},\]
we immediately obtain $\pa_{\underline{x}}\left[\fhi_{2n+1}^{m|0}\right]=\frac{\Gam\left(\frac{m}{2}-n\right)}{2^{2n} \, \pi^{m/2} \, \Gam(n)} \, \frac{|\underline{x}|^{2n}}{|\underline{x}|^m}= \fhi_{2n}^{m|0}$. This means that the above expression for $\fhi_{2n+1}^{m|0}$ constitutes a fundamental solution for $\pa_{\underline{x}}^{2n+1}$.

\noindent {Now,} substituting (\ref{PolyMonFundSolEven})-(\ref{PolyMonFundSolOdd}) into (\ref{FundSolSupDirOp})  we get
\begin{align}
\nu_1^{m|2n} &= \frac{\pi^n}{2\pi^{m/2}} \sum_{k=0}^{n-1}   \frac{\Gam\left(\frac{m}{2}-k-1\right)}{\Gam\left(n-k\right)}   \frac{|\underline{x}|^{2(k+1)}}{|\underline{x}|^m}\underline{x}\p^{\, 2n-2k-1} +  \frac{\pi^n}{2\pi^{m/2}} \sum_{k=0}^n  \frac{\Gam\left(\frac{m}{2}-k\right)}{\Gam\left(n-k+1\right)}   \frac{\underline{x} \, |\underline{x}|^{2k}}{|\underline{x}|^m}\underline{x}\p^{\, 2n-2k} \nonumber \\
&= \frac{1}{2\pi^{M/2}} \sum_{k=1}^{n}   \frac{\Gam\left(\frac{m}{2}-k\right)}{\Gam\left(n-k+1\right)}   \frac{|\underline{x}|^{2k}}{|\underline{x}|^m}\underline{x}\p^{\, 2n-2k+1} + \frac{1}{2\pi^{M/2}}  \sum_{k=0}^n  \frac{\Gam\left(\frac{m}{2}-k\right)}{\Gam\left(n-k+1\right)}   \frac{\underline{x} \, |\underline{x}|^{2k}}{|\underline{x}|^m}\underline{x}\p^{\, 2n-2k} \nonumber \\
&= \frac{1}{2\pi^{M/2}} \, {\bf x} \,\sum_{k=1}^n  \frac{\Gam\left(\frac{m}{2}-k\right)}{\Gam\left(n-k+1\right)} \,   |\underline{x}|^{2k-m} \, \underline{x}\p^{\, 2n-2k} + \frac{\Gam\left(\frac{m}{2}\right)}{2\pi^{M/2} \,\Gam\left(n+1\right)} \,  \frac{\underline{x}}{|\underline{x}|^m} \,\underline{x}\p^{\, 2n}. \label{IntFundSolFracForm}
\end{align}
{Recall} (see (\ref{PosPowBosElm})) that 
\begin{align*}
\frac{1}{|{\bf x}|^M} &= \left(|\underline{x}|^2-\underline{x}\p^{\, 2}\right)^{\frac{-M}{2}} = \sum_{j=0}^n \frac{\underline{x}\p^{\, 2j}}{j!} \, \frac{\Gam\left(\frac{m}{2}-m+j\right)}{\Gam\left(\frac{m}{2}-n \right)}  |\underline{x}|^{-m+2n-2j} \\
&= \frac{1}{\Gam\left(\frac{m}{2}-n \right)}  \sum_{k=0}^n  \frac{\Gam\left(\frac{m}{2}-k\right)}{\Gam(n-k+1)}  |\underline{x}|^{2k-m} \, \underline{x}\p^{\, 2n-2k}.
\end{align*}
Substituting {the later} into (\ref{IntFundSolFracForm}) we obtain 
\[\nu_1^{m|2n} = \frac{1}{2\pi^{M/2}} \, {\bf x} \left(\frac{\Gam\left(\frac{M}{2}\right)}{|{\bf x}|^M}- \frac{\Gam\left(\frac{m}{2}\right)}{\Gam\left(n+1\right)} \frac{\underline{x}\p^{\, 2n}}{|\underline{x}|^m} \right) + \frac{\Gam\left(\frac{m}{2}\right)}{2\pi^{M/2} \,\Gam\left(n+1\right)} \,  \frac{\underline{x}}{|\underline{x}|^m} \,\underline{x}\p^{\, 2n}= \frac{1}{|\Sa^{m-1|2n}|} \frac{{\bf x}}{|{\bf x}|^M},\]
which completes the proof.
$\hfill\square$

By doubling the bosonic dimension, i.e. $M=2m-2n$, we easily obtain the fundamental solution of $\pa_{{\bf J}({\bf x})}$ as shown in the next result.
\begin{cor}\label{FundSolHerJ}
For $M=2m-2n\notin -2\N+2$, i.e.\ $m>n$, the fundamental solutions of $\pa_{{\bf x}}$ and  $\pa_{{\bf J}({\bf x})}$ are given by,
\begin{align*}
\nu_1^{2m|2n} = \frac{1}{|\Sa^{2m-1|2n}|} \frac{{\bf x}}{|{\bf x}|^M},  \hspace{.3cm} \mbox{ and }   \hspace{.3cm} {\bf J}(\nu_1^{2m|2n}) = \frac{1}{|\Sa^{2m-1|2n}|} \frac{{\bf J}({\bf x})}{|{\bf x}|^M},
\end{align*}
respectively.
\end{cor}

\subsection{Distributional calculus for $\nu_1^{2m|2n}$ and ${\bf J}(\nu_1^{2m|2n})$}
One of the fundamental tools used in the distributional calculus with $\nu_1^{2m|2n}$ and ${\bf J}(\nu_1^{2m|2n})$ is the distribution {\it "finite part"} $\mbox{Fp}\, t^\mu_+$ on the real line. For a better understanding, we give its definition and list some of its main properties.

Let $\mu$ be a real parameter, $t$ a real variable and consider the function $t^\mu_+=\begin{cases} t^\mu, & t\geq0, \\ 0, & t<0. \end{cases}$ 
For $\mu>-1$ the function $t^\mu_+$ is locally integrable and hence constitutes a regular distribution, i.e.\ 
\begin{align}\label{RegDisttmu}
\langle t^\mu_+, \phi \rangle &= \int_{0}^{+\infty} t^\mu \phi(t)\, dt, & \phi&\in C^{\infty}_0(\R).
\end{align}
The finite part distribution $\mbox{Fp}\, t^\mu_+$ is an extension of the regular distribution $t^\mu_+$ to every value $\mu\in \R$. The idea of this extension is to consider only the finite part of the integral (\ref{RegDisttmu}). In this way we define
\[
\langle \mbox{Fp}\, t^\mu_+, \phi \rangle := \hspace{-.1cm}
\begin{cases}
\displaystyle \langle t^\mu_+, \phi \rangle, & -1<\mu, \\ \\

\displaystyle  \lim_{\ve\fd 0^+} \hspace{-.1cm} \left( \int_\ve^{+\infty} \hspace{-.1cm} t^{\mu} \phi(t)\, dt   + \sum_{j=1}^k \frac{\phi^{(j-1)}(0)}{(j-1)!}  \frac{\ve^{\mu+j}}{(\mu+j)}  \right), & \hspace{-.3cm} -(k+1) \hspace{-.1cm}<\mu < \hspace{-.1cm} -k, \\ \\

\displaystyle  \lim_{\ve\fd 0^+} \hspace{-.1cm}\left( \int_\ve^{+\infty} \hspace{-.1cm} t^{-k} \phi(t)\, dt   + \left(\sum_{j=1}^{k-1} \frac{\phi^{(j-1)}(0)}{(j-1)!}  \frac{\ve^{-k+j}}{(-k+j)}\right) + \frac{\phi^{(k-1)}(0)}{(k-1)!} \ln(\ve) \right), & \hspace{-.1cm} \mu=-k,
\end{cases}
\]
where $k\in \N$. 
The notation $ \mbox{\textup{Fp}}\, \int_{0}^{+\infty} t^\mu \phi(t)\, dt$ is often used for $\langle \mbox{Fp}\, t^\mu_+, \phi \rangle$.
\begin{pro}\label{FinPartProp}
 The following properties hold for $\mbox{\textup{Fp}}\, t^\mu_+$.
\begin{itemize}
\item[i)] $\displaystyle t \,  \mbox{\textup{Fp}}\, t^\mu_+ =  \mbox{\textup{Fp}}\, t^{\mu+1}_+$,  $ \mu\in\R$,
\item[ii)] $\displaystyle \frac{d}{dt} \mbox{\textup{Fp}}\, t^\mu_+ = \begin{cases} \displaystyle \mu  \mbox{\textup{Fp}}\, t^{\mu-1}_+,  & \mu\notin -\N+1, \\ \\[-3mm]
\displaystyle (-k)\, \mbox{\textup{Fp}}\,t^{-k-1}_+ + (-1)^k \frac{1}{k!} \del^{(k)}(t), & \mu=-k, \,k\in \N-1.\\ 
\end{cases}$
\end{itemize}
\end{pro}
In order to compute finite part distributions in $\R^m$ we need the so-called {\it generalized spherical means}, see e.g.\ \cite{MR2077080, MR1982060}. Let $\phi\in C^{\infty}_0(\R^m)$; putting $\underline{x}=r\underline{w}$, $r=|\underline{x}|$, we define the generalized spherical means
\begin{align*}
\Sigma^{(0)}[\phi](r) &= \frac{1}{|\Sa^{m-1}|} \int_{\Sa^{m-1}} \phi(r \underline{w}) \, dS_{\underline{w}}, &  \Sigma^{(1)}[\phi](r) &= \Sigma^{(0)}[\underline{w}\phi](r)=\frac{1}{|\Sa^{m-1}|} \int_{\Sa^{m-1}} \underline{w} \phi(r \underline{w}) \, dS_{\underline{w}}.
\end{align*}
It is clear that $\Sigma^{(0)}[\phi]: \R_+\fd \C$ and $\Sigma^{(1)}[\phi]: \R_+\fd \C^m$
are $C^\infty$ functions with singular support. 

We now list some important properties of these spherical means. The proofs of these results can be found in \cite{MR2077080}.
\begin{pro}\label{SpherMeansProp}
For a test function $\phi\in C^{\infty}_0(\R^m)$ one has
\begin{itemize}
\item[i)] $\Sigma^{(0)}[\underline{x}\, \phi]= r \, \Sigma^{(1)}[\phi]$,
\item[ii)] $\Sigma^{(1)}[\underline{x} \,\phi]= -r\, \Sigma^{(0)}[\phi]$,
\item[iii)] $\Sigma^{(0)}[\pa_{\underline{x}} \,\phi]= \left(\pa_r + \frac{m-1}{r}\right) \Sigma^{(1)}[\phi]$,
\item[iv)] $\Sigma^{(1)}[\pa_{\underline{x}} \,\phi]= -\pa_r \, \Sigma^{(0)}[\phi]$,
\item[v)] $\left\langle \del(r), \Sigma^{(0)}[\phi]\right\rangle= \langle \del(\underline{x}), \phi\rangle $,
\item[vi)] $\left\langle \del(r), \Sigma^{(1)}[\phi]\right\rangle=0$,
\item[vii)] $\left\langle \del' (r), \Sigma^{(1)}[\phi]\right\rangle=\frac{1}{m} \left\langle \pa_{\underline{x}} \del(\underline{x}) , \phi\right\rangle$.
\end{itemize}
\end{pro}

We now {have introduced all elements needed for} computing the action of the distribution $ \mbox{\textup{Fp}}\, |\underline{x}|^\lan_+$ on a test function $\phi\in C^{\infty}_0(\R^m)$, i.e.\
\begin{align}
\left\langle  \mbox{\textup{Fp}}\, |\underline{x}|^\lan_+, \phi\right\rangle &=  \mbox{\textup{Fp}}\, \int_{\R^m} |\underline{x}|^\lan \, \phi(\underline{x}) \, dV_{\underline{x}}=  \mbox{\textup{Fp}}\, \int_0^\infty \int_{\Sa^{m-1}} r^\lan \, \phi(r\underline{w}) r^{m-1} \, dr\, dS_{\underline{w}} \nonumber\\
&=   \mbox{\textup{Fp}}\, \int_0^\infty r^{\lan+m-1} \left( \int_{\Sa^{m-1}}  \phi(r\underline{w}) \, dS_{\underline{w}} \right)\, dr \nonumber \\
&= |\Sa^{m-1}| \; \mbox{\textup{Fp}}\, \int_0^\infty r^{\lan+m-1}  \Sigma^{(0)}[\phi](r) \, dr \nonumber \\
&= |\Sa^{m-1}| \left\langle  \mbox{\textup{Fp}}\,  r^{\lan+m-1}_+, \Sigma^{(0)}[\phi] \right\rangle. \label{FinPartinRm}
\end{align}
This motivates the introduction of the following distributions (see \cite{MR2077080, MR1982060}):
\begin{align*}
\left\langle  T_\lan, \phi \right\rangle &= |\Sa^{m-1}| \left\langle  \mbox{\textup{Fp}}\,  r^{\mu}_+ \, , \,\Sigma^{(0)}[\phi] \right\rangle, &
\left\langle  U_\lan, \phi \right\rangle &= |\Sa^{m-1}| \left\langle  \mbox{\textup{Fp}}\,  r^{\mu}_+ \,, \, \Sigma^{(1)}[\phi] \right\rangle,
\end{align*}
where $\mu= \lan+m-1$. In this way, one has in $\R^m$ that
\begin{align}\label{TUinTermsofFp}
 \mbox{\textup{Fp}}\, \frac{ |\underline{x}|^{\mu+1}}{ |\underline{x}|^m} &=T_\lan, &  \mbox{\textup{Fp}}\, \frac{ \underline{x} \, |\underline{x}|^{\mu}}{ |\underline{x}|^m} &=U_\lan.
\end{align}
Indeed, the first equality directly follows from (\ref{FinPartinRm}) while for the second one it suffices to note that
\begin{align*}
\left\langle  \mbox{\textup{Fp}}\, \frac{ \underline{x} \, |\underline{x}|^{\mu}}{ |\underline{x}|^m} , \phi\right\rangle &=  \mbox{\textup{Fp}}\, \int_{\R^m} \frac{ \underline{x} \, |\underline{x}|^{\mu}}{ |\underline{x}|^m}  \, \phi(\underline{x}) \, dV_{\underline{x}} = \mbox{\textup{Fp}}\, \int_0^\infty \int_{\Sa^{m-1}} r^{\mu+1-m} \, \underline{w} \,\phi(r\underline{w}) r^{m-1} \, dr\, dS_{\underline{w}} \\
&=  \mbox{\textup{Fp}}\, \int_0^\infty r^{\mu} \left( \int_{\Sa^{m-1}}  \underline{w} \, \phi(r\underline{w}) \, dS_{\underline{w}} \right)\, dr  \\
&= |\Sa^{m-1}| \left\langle  \mbox{\textup{Fp}}\,  r^{\mu}_+, \Sigma^{(1)}[\phi] \right\rangle \\
&= \left\langle  U_\lan, \phi \right\rangle.
\end{align*}
We can now compute  $\pa_{\bf x}[{\bf J}(\nu_1^{2m|2n})]$  and  $\pa_{{\bf J}({\bf x})}[\nu_1^{2m|2n}]$ in terms of the above defined distributions. 
\begin{pro}
For $m>n$ (i.e.\ $M=2(m-n)\notin -2\N+2$) it holds that
\begin{align}
\pa_{\bf x}[{\bf J}(\nu_1^{2m|2n})] &= \frac{2{\bf B}}{|\Sa^{2m-1|2n}|} \, \mbox{\textup{Fp}}\, \frac{1}{|{\bf x}|^M} + \frac{M}{|\Sa^{2m-1|2n}|} \, \mbox{\textup{Fp}}\,  \frac{{\bf x} {\bf J}({\bf x})}{{|{\bf x}|^{M+2}}} +  \frac{{\bf B}_b}{m} \del({\bf x}), \label{DerSupFundSol1}\\
 \pa_{{\bf J}({\bf x})}[\nu_1^{2m|2n}] &= \frac{-2{\bf B}}{|\Sa^{2m-1|2n}|} \, \mbox{\textup{Fp}}\, \frac{1}{|{\bf x}|^M} + \frac{M}{|\Sa^{2m-1|2n}|} \, \mbox{\textup{Fp}}\,  \frac{{\bf J}({\bf x}) {\bf x} }{{|{\bf x}|^{M+2}}} -  \frac{{\bf B}_b}{m} \del({\bf x}), \label{DerSupFundSol2}
\end{align}
where the distribution $\displaystyle \mbox{\textup{Fp}}\, \frac{1}{|{\bf x}|^{M+\al}}$, $\al\geq 0$, is defined in superspace as
\begin{align}\label{FpDistSS}
 \mbox{\textup{Fp}}\, \frac{1}{|{\bf x}|^{M+\al}}&= \frac{1}{\Gam(m-n+\frac{\al}{2})} \sum_{k=0}^n \frac{\Gam(m-k+\frac{\al}{2})}{\Gam(n-k+1)} \,  \mbox{\textup{Fp}}\, \frac{ |\underline{x}|^{2k}}{ |\underline{x}|^{2m+\al}} \, \underline{x}\p^{\, 2n-2k}.
\end{align}
\end{pro}
\begin{remark}
For n=0, formulae (\ref{DerSupFundSol1})-(\ref{DerSupFundSol2}) coincide with the {corresponding} expressions in the purely bosonic case computed in \cite{MR2540516}.
\end{remark}
\pf
From Corollary \ref{FundSolHerJ} we have in distributional sense that 
\[{\bf J}(\nu_1^{2m|2n})= \frac{1}{|\Sa^{2m-1|2n}|} {\bf J}({\bf x}) \, \mbox{\textup{Fp}}\,  \frac{1}{|{\bf x}|^{M}} . \]
Hence, 
\begin{align}
\pa_{\bf x}[{\bf J}(\nu_1^{2m|2n})]  &=  \frac{1}{|\Sa^{2m-1|2n}|} \left( \pa_{\bf x}\left[ \mbox{\textup{Fp}}\,  \frac{1}{|{\bf x}|^{M}}\right] {\bf J}({\bf x}) \, + \, \mbox{\textup{Fp}}\,  \frac{1}{|{\bf x}|^{M}} \pa_{\bf x} \left[{\bf J}({\bf x})\right] \right) \nonumber \\
&= \frac{2{\bf B}}{|\Sa^{2m-1|2n}|} \, \mbox{\textup{Fp}}\, \frac{1}{|{\bf x}|^M} \, + \,  \frac{1}{|\Sa^{2m-1|2n}|} \, \pa_{\bf x}\left[ \mbox{\textup{Fp}}\,  \frac{1}{|{\bf x}|^{M}}\right] {\bf J}({\bf x}). \label{FirstGenExp}
\end{align}
Using (\ref{TUinTermsofFp}) and (\ref{FpDistSS}) we obtain,
\begin{equation}\label{DetFpBoldXM}
\pa_{\bf x}\left[ \mbox{\textup{Fp}}\,  \frac{1}{|{\bf x}|^{M}}\right] =  \frac{1}{\Gam(m-n)}  \sum_{k=0}^n \frac{\Gam(m-k)}{\Gam(n-k+1)} \, \pa_{\bf x}\left[  T_{2k-2m}\, \underline{x}\p^{\, 2n-2k} \right],
\end{equation}
where
\begin{align*}
 \pa_{\bf x}\left[  T_{2k-2m}\, \underline{x}\p^{\, 2n-2k} \right] =  -\pa_{\underline{x}}\left[  T_{2k-2m} \right] \, \underline{x}\p^{\, 2n-2k} +  T_{2k-2m}  \, \pa_{\underline{x}\p}\left[ \underline{x}\p^{\, 2n-2k} \right].
\end{align*}
By Propositions \ref{FinPartProp} and \ref{SpherMeansProp} we now get,
\begin{align*}
\langle  \pa_{\underline{x}}\left[  T_{2k-2m} \right], \phi\rangle &= -\langle  T_{2k-2m} ,   \pa_{\underline{x}}\left[ \phi\right] \rangle = -|\Sa^{2m-1}| \left\langle  \mbox{\textup{Fp}}\,  r^{2k-1}_+ \, , \,\Sigma^{(0)}\left[  \pa_{\underline{x}}\phi\right] \right\rangle \\
&= -|\Sa^{2m-1}| \left\langle  \mbox{\textup{Fp}}\,  r^{2k-1}_+ \, , \,\left(\pa_r + \frac{2m-1}{r}\right) \Sigma^{(1)}[\phi]\right\rangle \\
&= |\Sa^{2m-1}| \left\langle \frac{d}{dr} \, \mbox{\textup{Fp}}\,  r^{2k-1}_+ \, ,\Sigma^{(1)}[\phi]\right\rangle -|\Sa^{2m-1}| \left\langle (2m-1)  \, \mbox{\textup{Fp}}\,  r^{2k-2}_+ \, ,  \Sigma^{(1)}[\phi]\right\rangle \\
&= |\Sa^{2m-1}| \left\langle (2k-1)  \, \mbox{\textup{Fp}}\,  r^{2k-2}_+ \, -\del_{k,0} \del'(r) - (2m-1)  \, \mbox{\textup{Fp}}\,  r^{2k-2}_+ \,,  \Sigma^{(1)}[\phi]\right\rangle \\
&=(2k-2m) |\Sa^{2m-1}|  \left\langle  \mbox{\textup{Fp}}\,  r^{2k-2}_+ \, , \, \Sigma^{(1)}[\phi]\right\rangle - \del_{k,0} |\Sa^{2m-1}|  \left\langle   \del'(r) \, , \, \Sigma^{(1)}[\phi]\right\rangle \\
&=   \left\langle (2k-2m)U_{2k-2m-1} - \frac{ \del_{k,0} |\Sa^{2m-1}|}{2m} \pa_{\underline{x}} \del(\underline{x})    \, , \, \phi \right\rangle, 
\end{align*}
or equivalently,
\begin{equation*}
\pa_{\underline{x}}\left[  T_{2k-2m} \right] = (2k-2m)U_{2k-2m-1} - \frac{ \del_{k,0} |\Sa^{2m-1}|}{2m} \pa_{\underline{x}} \del(\underline{x}).  
\end{equation*}
Moreover,
\begin{align*}
 \pa_{\underline{x}\p}\left[ \underline{x}\p^{\, 2n-2k} \right]&= 2(n-k)  \underline{x}\p^{\, 2n-2k-1}, & k&=0, 1,\ldots, n,
\end{align*}
where we are formally defining $0 \underline{x}\p^{\,-1} :=0$ in the case $k=n$ (we recall that the element  $\underline{x}\p^{\,-1}$  does not exist due to the nilpotency of $ \underline{x}\p$}). Hence we conclude that
\[\pa_{\bf x}\left[  T_{2k-2m}\, \underline{x}\p^{\, 2n-2k} \right] = 2(m-k)U_{2k-2m-1} \,  \underline{x}\p^{\, 2n-2k}  + 2(n-k) \, T_{2k-2m} \, \underline{x}\p^{\, 2n-2k-1} + \frac{ \del_{k,0} |\Sa^{2m-1}|}{2m} \pa_{\underline{x}} \del(\underline{x}) \,  \underline{x}\p^{\, 2n} .
\]
Substituting this into (\ref{DetFpBoldXM}) we obtain,
\begin{align*}
\pa_{\bf x}\left[ \mbox{\textup{Fp}}\,  \frac{1}{|{\bf x}|^{M}}\right] &=  \frac{1}{\Gam(m-n)}  \sum_{k=0}^n \frac{2\Gam(m-k)}{\Gam(n-k+1)} \, \left[(m-k)U_{2k-2m-1} \,  \underline{x}\p^{\, 2n-2k} + (n-k) \, T_{2k-2m} \, \underline{x}\p^{\, 2n-2k-1}\right] \\
&\phantom{=} + \frac{\Gam(m)}{\Gam(m-n) \Gam(n+1)} \frac{|\Sa^{2m-1}|}{2m} \pa_{\underline{x}} \del(\underline{x}) \,  \underline{x}\p^{\, 2n} \\
&=  \frac{2}{\Gam(m-n)} \left[ \sum_{k=0}^n \frac{\Gam(m-k+1)}{\Gam(n-k+1)} \, U_{2k-2m-1}  \,  \underline{x}\p^{\, 2n-2k} + \sum_{k=1}^{n} \frac{\Gam(m-k+1)}{\Gam(n-k+1)} \, T_{2k-2-2m} \, \underline{x}\p^{\, 2n-2k+1} \right] \\
&\phantom{=} +\frac{\pi^m}{m \Gam(m-n) \Gam(n+1)}\, \pa_{\underline{x}} \del(\underline{x}) \,  \underline{x}\p^{\, 2n} \\
&=  \frac{2}{\Gam(m-n)}  \sum_{k=1}^n \frac{\Gam(m-k+1)}{\Gam(n-k+1)} \left( \mbox{\textup{Fp}}\, \frac{\underline{x} |\underline{x}|^{2k-2}}{ |\underline{x}|^{2m}} +  \mbox{\textup{Fp}}\, \frac{\underline{x}\p |\underline{x}|^{2k-2}}{|\underline{x}|^{2m}}\right) \underline{x}\p^{\, 2n-2k} \\
&\phantom{=} +  \frac{2\Gam(m+1)}{\Gam(m-n) \Gam(n+1)} U_{-2m-1}  \,  \underline{x}\p^{\, 2n} +\frac{\pi^m}{m \Gam(m-n) \Gam(n+1)}\, \pa_{\underline{x}} \del(\underline{x}) \,  \underline{x}\p^{\, 2n} \\
&=   \frac{2 {\bf x}}{\Gam(m-n)}  \sum_{k=1}^n \left(\frac{\Gam(m-k+1)}{\Gam(n-k+1)} \, \mbox{\textup{Fp}}\, \frac{|\underline{x}|^{2k}}{ |\underline{x}|^{2m+2}} \, \underline{x}\p^{\, 2n-2k} \right)+  \frac{2\Gam(m+1)}{\Gam(m-n) \Gam(n+1)} U_{-2m-1}  \,  \underline{x}\p^{\, 2n} \\
&\phantom{=} +\frac{\pi^m}{m \Gam(m-n) \Gam(n+1)}\, \pa_{\underline{x}} \del(\underline{x}) \,  \underline{x}\p^{\, 2n} .
\end{align*}
Then from (\ref{FpDistSS}) we get,
\begin{align}
\pa_{\bf x}\left[ \mbox{\textup{Fp}}\,  \frac{1}{|{\bf x}|^{M}}\right] &=  \frac{2 {\bf x}}{\Gam(m-n)}  \left[ \Gam(m-n+1) \mbox{\textup{Fp}}\, \frac{1}{|{\bf x}|^{M+2}} -\frac{\Gam(m+1)}{\Gam(n+1)} \,  \mbox{\textup{Fp}}\, \frac{ 1}{ |\underline{x}|^{2m+2}} \, \underline{x}\p^{\, 2n} \right] \nonumber \\
&\phantom{=} +  \frac{2\Gam(m+1)}{\Gam(m-n) \Gam(n+1)} U_{-2m-1}  \,  \underline{x}\p^{\, 2n} +\frac{\pi^m}{m \Gam(m-n) \Gam(n+1)}\, \pa_{\underline{x}} \del(\underline{x}) \,  \underline{x}\p^{\, 2n} \nonumber\\
&= 2(m-n) \, \mbox{\textup{Fp}}\, \frac{{\bf x}}{|{\bf x}|^{M+2}} +\frac{\pi^m}{m \Gam(m-n) \Gam(n+1)}\, \pa_{\underline{x}} \del(\underline{x}) \,  \underline{x}\p^{\, 2n}. \label{DerFpModBfXM}
\end{align}
{Substitution of} (\ref{DerFpModBfXM}) into (\ref{FirstGenExp}) {yields}
\begin{align}\label{AlmFinComforDerJnu}
\pa_{\bf x}[{\bf J}(\nu_1^{2m|2n})]  &=  \frac{1}{|\Sa^{2m-1|2n}|} \left( 2{\bf B} \, \mbox{\textup{Fp}}\, \frac{1}{|{\bf x}|^M} \, + \,  M\, \mbox{\textup{Fp}}\, \frac{{\bf x} {\bf J}({\bf x})}{|{\bf x}|^{M+2}} +\frac{\pi^m}{m \Gam(m-n) \Gam(n+1)}\, \pa_{\underline{x}} \del(\underline{x}) \,  {\bf J}(\underline{x}) \,\underline{x}\p^{\, 2n} \right).
\end{align}
Observe now that, for any complex valued test function $\phi$, one has
\begin{align*}
\left \langle   \pa_{\underline{x}} \del(\underline{x}) \,  {\bf J}(\underline{x}) , \phi \right\rangle &= -\left \langle    \del(\underline{x}) \,  ,  \pa_{\underline{x}} \left[{\bf J}(\underline{x})\phi\right] \right\rangle = - \left \langle    \del(\underline{x}) \,  ,  -2{\bf B}_b \phi +  \, \pa_{\underline{x}} \left[\phi\right] {\bf J}(\underline{x})\right\rangle =  \left \langle   2{\bf B}_b \del(\underline{x}) \,  ,   \phi \right\rangle,
\end{align*}
which implies that  $\pa_{\underline{x}} \del(\underline{x}) \,  {\bf J}(\underline{x}) = 2{\bf B}_b \del(\underline{x})$. Hence, since $\del({\bf x})= \del(\underline{x}) \frac{\pi^n}{n!}\underline{x}\p^{\, 2n}$, we obtain,
\begin{align*}
\frac{\pi^m}{m \Gam(m-n) \Gam(n+1)}\, \pa_{\underline{x}} \del(\underline{x}) \,  {\bf J}(\underline{x}) \,\underline{x}\p^{\, 2n}= \frac{2\pi^m}{m \Gam(\frac{M}{2}) n!}\, {\bf B}_b \del(\underline{x})\underline{x}\p^{\, 2n} = \frac{|\Sa^{2m-1|2n}|}{m}  {\bf B}_b \del({\bf x}).
\end{align*}
Finally, substituting {the latter} into (\ref{AlmFinComforDerJnu}) we obtain (\ref{DerSupFundSol1}). Formula (\ref{DerSupFundSol2}) easily follows from applying ${\bf J}$ to both sides of (\ref{DerSupFundSol1}) and using the properties ${\bf J}^2[\nu_1^{2m|2n}]=-\nu_1^{2m|2n}$ and ${\bf J}({\bf B})={\bf B}$.
$\hfill\square$

\subsection{Hermitian counterparts of $\nu_1^{2m|2n}$ and ${\bf J}(\nu_1^{2m|2n})$}
Similarly as above, we introduce the following Hermitian counterparts to the pair of fundamental solutions $(\nu_1^{2m|2n},{\bf J}(\nu_1^{2m|2n}))$, for $m>n$:
\begin{align*}
\Psi^{m|n}_1 &= \nu_1^{2m|2n}+ i {\bf J}(\nu_1^{2m|2n}), & {\Psi^{m|n}_1}^\dag &= -\left(\nu_1^{2m|2n}- i {\bf J}(\nu_1^{2m|2n}) \right),
\end{align*}
or equivalently,
\begin{align*}
\Psi^{m|n}_1({\bf Z}) &= \frac{2}{|\Sa^{2m-1|2n}|} \frac{{\bf Z}}{|{\bf Z}|^M}, & {\Psi^{m|n}_1}^\dag ({\bf Z})&= \frac{2}{|\Sa^{2m-1|2n}|} \frac{{\bf Z}^\dag}{|{\bf Z}|^M},
\end{align*}
where we recall that $|{\bf Z}|=|{\bf x}|$. As in the purely bosonic case, see e.g. \cite{MR2540516}, $\Psi^{m|n}_1$ and $ {\Psi^{m|n}_1}^\dag$ are not the fundamental solutions of the Hermitian super Dirac operators $\pa_{\bf Z}$ and $\pa_{{\bf Z}^\dag}$. Indeed, from (\ref{DerSupFundSol1})-(\ref{DerSupFundSol2}) one obtains the following results. 
\begin{pro}\label{DerHermCountAll}
\begin{align}
\pa_{\bf Z} \Psi^{m|n}_1 &= \frac{m+i{\bf B}_b}{2m} \del({\bf x}) + \frac{\frac{M}{2} +i{\bf B}}{|\Sa^{2m-1|2n}|}  \, \mbox{\textup{Fp}}\, \frac{1}{|{\bf x}|^M} - \frac{M}{|\Sa^{2m-1|2n}|} \, \mbox{\textup{Fp}}\, \frac{{\bf Z}^\dag {\bf Z}}{|{\bf x}|^{M+2}} =  {\Psi^{m|n}_1}^\dag  \pa_{{\bf Z}^\dag} , \label{DerHermCout1} \\
\pa_{{\bf Z}^\dag} \Psi^{m|n}_1 &=0=  {\Psi^{m|n}_1}^\dag  \pa_{{\bf Z}} ,  \label{DerHermCout2} \\
 \pa_{{\bf Z}} {\Psi^{m|n}_1}^\dag &= 0  =  \Psi^{m|n}_1 \pa_{{\bf Z}^\dag}, \nonumber \\
 \pa_{{\bf Z}^\dag} {\Psi^{m|n}_1}^\dag &= \frac{m-i{\bf B}_b}{2m} \del({\bf x}) - \frac{\frac{M}{2} +i{\bf B}}{|\Sa^{2m-1|2n}|}  \, \mbox{\textup{Fp}}\, \frac{1}{|{\bf x}|^M} + \frac{M}{|\Sa^{2m-1|2n}|} \, \mbox{\textup{Fp}}\, \frac{{\bf Z}^\dag {\bf Z}}{|{\bf x}|^{M+2}} =  {\Psi^{m|n}_1}  \pa_{{\bf Z}}. \nonumber
\end{align}
\end{pro}
\pf
We will only prove the left equalities in (\ref{DerHermCout1}) and (\ref{DerHermCout2}) since the {remaining ones} can be proven along similar lines.  We first observe that 
\begin{align*}
\pa_{\bf Z} \Psi^{m|n}_1 &= \frac{1}{4}\left(\pa_{\bf x}-i \pa_{\bf J({\bf x})}\right)\left( \nu_1^{2m|2n}+ i {\bf J}(\nu_1^{2m|2n})\right)\\
&=  \frac{1}{4} \left[ \left(\pa_{\bf x} \nu_1^{2m|2n} + \pa_{\bf J({\bf x})}{\bf J}(\nu_1^{2m|2n})\right) + i \left( \pa_{\bf x}{\bf J}(\nu_1^{2m|2n})   - \pa_{\bf J({\bf x})} \nu_1^{2m|2n} \right)  \right] \\
&= \frac{1}{4} \left[ 2\del({\bf x}) +i \left( \frac{4{\bf B}}{|\Sa^{2m-1|2n}|} \, \mbox{\textup{Fp}}\, \frac{1}{|{\bf x}|^M} + \frac{2M}{|\Sa^{2m-1|2n}|} \, \mbox{\textup{Fp}}\,  \frac{{\bf x} {\bf J}({\bf x})}{{|{\bf x}|^{M+2}}} +  \frac{2{\bf B}_b}{m} \del({\bf x}) \right) \right] \\
&=  \frac{m+i{\bf B}_b}{2m} \del({\bf x}) + \frac{i{\bf B}}{|\Sa^{2m-1|2n}|} \, \mbox{\textup{Fp}}\, \frac{1}{|{\bf x}|^M}  + \frac{iM}{2|\Sa^{2m-1|2n}|} \, \mbox{\textup{Fp}}\,  \frac{{\bf x} {\bf J}({\bf x})}{{|{\bf x}|^{M+2}}}.
\end{align*}
On the other hand it can be easily proven that ${\bf x} {\bf J}({\bf x})=-i|{\bf x}|^2 +2i{\bf Z}^\dag {\bf Z}$. Substituting this {result} into the above formula, we get (\ref{DerHermCout1}). For (\ref{DerHermCout2}) it suffices to note that
\begin{align*}
\pa_{{\bf Z}^\dag} \Psi^{m|n}_1 &=  -\frac{1}{4}\left(\pa_{\bf x}+i \pa_{\bf J({\bf x})}\right)\left( \nu_1^{2m|2n}+ i {\bf J}(\nu_1^{2m|2n})\right)\\
& = - \frac{1}{4} \left[ \left(\pa_{\bf x} \nu_1^{2m|2n} - \pa_{\bf J({\bf x})}{\bf J}(\nu_1^{2m|2n})\right) + i \left( \pa_{\bf x}{\bf J}(\nu_1^{2m|2n})   + \pa_{\bf J({\bf x})} \nu_1^{2m|2n} \right)  \right] \\
&=0,
\end{align*}
which completes the proof.
$\hfill\square$

Proposition \ref{DerHermCountAll} shows that the functions $ \Psi^{m|n}_1$ and $ {\Psi^{m|n}_1}^\dag $ are not sh-monogenic. Nevertheless, they can be combined in a ($2\times 2$) circulant matrix in order to obtain the Hermitian Cauchy formulae in superspace. This {approach is} inspired {by the one} used in the purely bosonic case, see e.g.\ \cite{MR2540516, MR1889406}.
\begin{teo}\label{FisrtMatApp}
Introducing the particular circulant $\textup{(}2\times 2\textup{)}$ matrices 
\begin{align*}
{\bm{\mathcal D}}_{({\bf Z}, {\bf Z}^\dag)} &= \left(\begin{array}{cc} \pa_{\bf Z} & \pa_{{\bf Z}^\dag} \\ \pa_{{\bf Z}^\dag} & \pa_{{\bf Z}}\end{array}\right), 
&
{\bm \Psi}^{m|n}_{2\times 2} &= \left(\begin{array}{cc}{\Psi^{m|n}_1} &{\Psi^{m|n}_1}^\dag \\ {\Psi^{m|n}_1}^\dag & {\Psi^{m|n}_1} \end{array}\right),
&
{\bm \del} &:= \del I_2= \left(\begin{array}{cc} \del & 0 \\0 & \del \end{array}\right),
\end{align*}
one obtains that 
\[{\bm{\mathcal D}}_{({\bf Z}, {\bf Z}^\dag)}\; {\bm \Psi}^{m|n}_{2\times 2}({\bf Z}) = {\bm \del}({\bf x}) =  {\bm \Psi}^{m|n}_{2\times 2}({\bf Z}) \; {\bm{\mathcal D}}_{({\bf Z}, {\bf Z}^\dag)}.\]
Here, $I_2$ denotes the identity matrix of order $\textup{(}2\times 2\textup{)}$.
\end{teo}

\section{Hermitian Cauchy-Pompeiu formula in superspace}
Theorem \ref{FisrtMatApp} means that ${\bm \Psi}^{m|n}_{2\times 2}$ may be considered as a fundamental solution of ${\bm{\mathcal D}}_{({\bf Z}, {\bf Z}^\dag)}$ in the above-introduced matrix context. This observation is crucial for the matrix approach used in Hermitian Clifford analysis to arrive at a Cauchy integral formula. Moreover, it is remarkable that the Dirac matrix ${\bm{\mathcal D}}_{({\bf Z}, {\bf Z}^\dag)}$ factorizes in some sense the Laplacian, i.e.\
\begin{align*}
{\bm{\mathcal D}}_{({\bf Z}, {\bf Z}^\dag)} \; \left({\bm{\mathcal D}}_{({\bf Z}, {\bf Z}^\dag)}\right)^\dag &= \frac{1}{4} \left(\begin{array}{cc} \Del_{2m|2n} & 0 \\ 0 & \Del_{2m|2n} \end{array}\right), & \mbox{with } && \left({\bm{\mathcal D}}_{({\bf Z}, {\bf Z}^\dag)}\right)^\dag & =  \left(\begin{array}{cc} \pa_{{\bf Z}^\dag} & \pa_{{\bf Z}} \\ \pa_{{\bf Z}} & \pa_{{\bf Z}^\dag} \end{array}\right).
\end{align*}
Thus, in the same setting, we associate, with every pair of Clifford-valued superfunctions $G_1, G_2\in C^\infty(\Om)\otimes \mathfrak{G}_{2n}\otimes \mathcal{C}_{2m,2n}$, the matrix function 
\begin{equation}\label{CircMatGenForm}
{\bm G}_2^1=  \left(\begin{array}{cc} G_1& G_2 \\ G_2 & G_1 \end{array}\right).
\end{equation}   
\begin{defi}
The matrix function ${\bm G}_2^1$ is said to be (left) ${\bm{SH}}$-monogenic if ${\bm{\mathcal D}}_{({\bf Z}, {\bf Z}^\dag)} {\bm G}_2^1 = \bm{0}$, where $\bm{0}$ denotes the matrix with zero entries. 
\end{defi}
The above definition for  ${\bm{SH}}$-monogenicity explicitly reads 
\[
\begin{cases}
\pa_{\bf Z} [G_1]+ \pa_{{\bf Z}^\dag} [G_2] =0,\\
\pa_{\bf Z} [G_2]+ \pa_{{\bf Z}^\dag} [G_1] =0.
\end{cases}
\]
When considering in particular  $G_1=G$ and $G_2=G^\dag$, the ${\bm{SH}}$-monogenicity of the corresponding matrix function ${\bm G}=\left(\begin{array}{cc} G& G^\dag \\ G^\dag & G \end{array}\right)$ does not imply, in general, the sh-monogenicity of $G$ and vice versa. As a clear example to illustrate this consider the matrix ${\bm G}={\bm \Psi}^{m|n}_{2\times 2}$, i.e.\ $G= \Psi^{m|n}_1$. An important exception to this general remark occurs in the case of Grassmann-valued functions. Indeed,  if $G\in C^\infty (\Om) \otimes \mathfrak{G}_{2n}$ one has
\begin{align*} 
\pa_{\bf Z} [G]+ \pa_{{\bf Z}^\dag} [G^\dag] &= \sum_{j=1}^m \left(\f_j^\dag \,\pa_{z_j}[G]+ \f_j \,\pa_{z_j^c}[G^c] \right)+ 2i\sum_{j=1}^n \left(\f_j\p^\dag \, \pa_{z\p_j}[G]-\f_j\p \, \pa_{z\p_j^c}[G^c]\right), \\
\pa_{\bf Z} [G^\dag]+ \pa_{{\bf Z}^\dag} [G] &= \sum_{j=1}^m \left(\f_j^\dag \,\pa_{z_j}[G^c]+ \f_j \,\pa_{z_j^c}[G] \right)+ 2i\sum_{j=1}^n \left(\f_j\p^\dag \, \pa_{z\p_j}[G^c]-\f_j\p \, \pa_{z\p_j^c}[G]\right),
\end{align*}
which means in this case that 
\[
\begin{cases}
\pa_{\bf Z} [G]+ \pa_{{\bf Z}^\dag} [G^\dag] =0,\\
\pa_{\bf Z} [G^\dag]+ \pa_{{\bf Z}^\dag} [G] =0,
\end{cases}
\iff
\begin{cases}
\pa_{z_j}[G]=  \pa_{z_j^c}[G]=0, & j=1,\ldots, m,\\
\pa_{z\p_j}[G] =\pa_{z\p_j^c}[G]= 0, & j=1,\ldots, n,
\end{cases}
\iff
\pa_{\bf Z} [G]= 0=\pa_{{\bf Z}^\dag} [G]. 
\]
Another important case occurs when considering the matrix function ${\bm G}_0= GI_2=\left(\begin{array}{cc} G& 0 \\ 0 & G \end{array}\right)$ with $G\in C^\infty(\Om)\otimes \mathfrak{G}_{2n}\otimes \mathcal{C}_{2m,2n}$. {Also} in this case the ${\bm{SH}}$-monogenicity of ${\bm G}_0$ is equivalent to the sh-monogenicity of $G$. Indeed,
\[{\bm{\mathcal D}}_{({\bf Z}, {\bf Z}^\dag)}{\bm G}_0 = \left(\begin{array}{cc} \pa_{{\bf Z}}[G]& \pa_{{\bf Z}^\dag}[G]\\ \pa_{{\bf Z}^\dag}[G] & \pa_{{\bf Z}}[G] \end{array}\right).\]
We are now in conditions {to reformulate} the Hermitian-Stokes theorem given in Corollary \ref{Her_stokThemCor2}, in a matrix form. The proof easily follows by taking deliberate combinations of the formulae (\ref{Stok11})-(\ref{Stok22}).
\begin{teo}
Let ${\bm F}_2^1$ and ${\bm G}_2^1$ be a pair of  matrix functions  of the form (\ref{CircMatGenForm}) with entries in $C^\infty(\Om)\otimes \mathfrak{G}_{2n}\otimes \mathcal{C}_{2m,2n}$. Let moreover $\al$ and $\be$ be distributions in $ \mathcal{E}' \otimes \mathfrak{G}_{2n}^{(ev)}$ and consider the circulant matrix distribution ${\bm \Sigma}=\left(\begin{array}{cc} \al& \be \\ \be & \al \end{array}\right)$ such that $supp \,\Sigma:=supp \, \al\cup supp \, \be$
 is a subset of $\Om$.  It then holds that
\begin{equation}\label{MatHerStokSupSDist}
\int_{\R^{2m|2n}} \left({\bm F}_2^1{\bm{\mathcal D}}_{({\bf Z}, {\bf Z}^\dag)} \right) \, {\bm \Sigma} \, {\bm G}_2^1 \,+\, {\bm F}_2^1 \, {\bm \Sigma} \, \left( {\bm{\mathcal D}}_{({\bf Z}, {\bf Z}^\dag)} {\bm G}_2^1 \right) = -\int_{\R^{2m|2n}} {\bm F}_2^1 \, \left({\bm{\mathcal D}}_{({\bf Z}, {\bf Z}^\dag)} {\bm \Sigma}\right) \, {\bm G}_2^1.
\end{equation}
\end{teo}
{Considering now} a phase function $g=g_0+{\bf g}\in C^\infty(\R^{2m})\otimes \mathfrak{G}_{2n}^{(ev)}$ such that $\{g_0\leq 0\}\inc\Om$ is compact and $\pa_{\underline x}[g_0]\neq0$ on $g_0^{-1}(0)$, the corresponding matrix distribution 
\begin{equation}\label{MatSubsStok}
 {\bm \Sigma}= H(-g)I_2= \left(\begin{array}{cc} H(-g)& 0 \\ 0 & H(-g) \end{array}\right)
 \end{equation}
is a commuting matrix {satisfying} 
\[{\bm{\mathcal D}}_{({\bf Z}, {\bf Z}^\dag)} {\bm \Sigma} = -\del(g) \left(\begin{array}{cc} \pa_{{\bf Z}}[g]& \pa_{{\bf Z}^\dag}[g]\\ \pa_{{\bf Z}^\dag}[g] & \pa_{{\bf Z}}[g] \end{array}\right)=-\del(g)  {\bm{\mathcal D}}_{({\bf Z}, {\bf Z}^\dag)}[g]. \]
Hence formula (\ref{MatHerStokSupSDist}) takes the following form, when {substituting} (\ref{MatSubsStok}):
\begin{equation}\label{StokHeavmatSupSHer}
\int_{\R^{2m|2n}} H(-g) \left[ \left({\bm F}_2^1{\bm{\mathcal D}}_{({\bf Z}, {\bf Z}^\dag)} \right) \,  {\bm G}_2^1 \,+\, {\bm F}_2^1 \, \left( {\bm{\mathcal D}}_{({\bf Z}, {\bf Z}^\dag)} {\bm G}_2^1 \right)  \right] =  \int_{\R^{2m|2n}} {\bm F}_2^1 \,\del(g)  {\bm{\mathcal D}}_{({\bf Z}, {\bf Z}^\dag)}[g] \, {\bm G}_2^1.
\end{equation}
In order to prove the Hermitian Cauchy-Pompeiu formula in superspace we first observe that in (\ref{MatHerStokSupSDist}) the matrix function ${\bm F}_2^1$ can be replaced by any matrix distribution 
\begin{align}\label{MatDistGenForm}
{\bm \Upsilon}&=\left(\begin{array}{cc} \gam& \sigma\\ \sigma & \gam \end{array}\right), & \gam, \sigma &\in \mathcal{E}'\otimes \mathfrak{G} _{2n}\otimes  \mathcal C_{2m,2n},
\end{align}
such that the sets
\begin{align*}
sing \; supp \, {\bm \Upsilon}&:= sing \; supp \, \gam \cup sing \; supp \, \sigma & \mbox{ and } && sing \; supp \, {\bm \Sigma}&:= sing \; supp \, \al \cup sing \; supp \, \be
\end{align*}
are disjoint. 
Under these conditions, (\ref{MatHerStokSupSDist}) can be proven for ${\bm F}_2^1 = {\bm \Upsilon}$ by  taking deliberate combinations of (\ref{Stok11})-(\ref{Stok22}) {with} $F=\gam$ and $F=\sigma$. These substitutions are possible since distributions  with disjoint singular supports can be multiplied and the Leibniz rule {remains valid} for such a product, see (\ref{DistProd}) and (\ref{prodSupDist}). Applying this reasoning to (\ref{StokHeavmatSupSHer}), for which we are {taking} $\Sigma=H(-g)I_2$, we immediately obtain the following consequence.
\begin{cor}
Let ${\bm G}_2^1$ be a matrix function of the form (\ref{CircMatGenForm}) with entries in $C^\infty(\Om)\otimes \mathfrak{G}_{2n}\otimes \mathcal{C}_{2m,2n}$. Let $g=g_0+{\bf g}\in C^\infty(\R^{2m})\otimes \mathfrak{G}_{2n}^{(ev)}$ be a phase function such that $\{g_0\leq 0\}\inc \Om$ is compact and $\pa_{\underline x}[g_0]\neq0$ on $g_0^{-1}(0)$, and let ${\bm \Upsilon}$ be a matrix distribution of the form (\ref{MatDistGenForm}) such that $sing \; supp \, {\bm \Upsilon} \cap g^{-1}(0)=\emptyset$. It then holds that
\begin{equation}\label{StokHeavmatSupSHerDist}
\int_{\R^{2m|2n}} H(-g) \left[ \left({\bm {\Upsilon}}{\bm{\mathcal D}}_{({\bf Z}, {\bf Z}^\dag)} \right) \,  {\bm G}_2^1 \,+\,{\bm \Upsilon} \, \left( {\bm{\mathcal D}}_{({\bf Z}, {\bf Z}^\dag)} {\bm G}_2^1 \right)  \right] =  \int_{\R^{2m|2n}} {\bm \Upsilon} \,\del(g)  {\bm{\mathcal D}}_{({\bf Z}, {\bf Z}^\dag)}[g] \, {\bm G}_2^1,
\end{equation}
where the distributional products $H(-g)\left({\bm {\Upsilon}}{\bm{\mathcal D}}_{({\bf Z}, {\bf Z}^\dag)} \right)$, $H(-g){\bm \Upsilon}$ and ${\bm \Upsilon} \del(g)$ are {to be} understood in the sense of (\ref{DistProd}) and (\ref{prodSupDist}).
\end{cor}
\pf It suffices to note that $sing \; supp \, H(-g)= g^{-1}(0)$. $\hfill\square$

Let us now consider the supervector ${\bf y}=\underline{y}+\underline{y}\p$, its Hermitian counterparts 
\begin{align*}
{\bf U} &=\frac{1}{2}({\bf y}+i{\bf J}({\bf y})), & {\bf U}^\dag &=-\frac{1}{2}({\bf y}-i{\bf J}({\bf y})),
\end{align*}
and the matrix distribution 
\[{\bm \Psi}^{m|n}_{2\times 2} ({\bf Z}-{\bf U})= \left(\begin{array}{cc}{\Psi^{m|n}_1} ({\bf Z}-{\bf U}) &{\Psi^{m|n}_1}^\dag ({\bf Z}-{\bf U}) \\ {\Psi^{m|n}_1}^\dag ({\bf Z}-{\bf U}) & {\Psi^{m|n}_1} ({\bf Z}-{\bf U}) \end{array}\right),\]
where we recall that 
\begin{align*}
\Psi^{m|n}_1({\bf Z}-{\bf U}) &= \frac{2}{|\Sa^{2m-1|2n}|} \frac{{\bf Z}-{\bf U}}{|{\bf Z}-{\bf U}|^M}, & {\Psi^{m|n}_1}^\dag ({\bf Z}-{\bf U})&= \frac{2}{|\Sa^{2m-1|2n}|} \frac{{\bf Z}^\dag-{\bf U}^\dag}{|{\bf Z}-{\bf U}|^M}.
\end{align*}
The following Hermitian Cauchy-Pompeiu formula in superspace {then} is established.
\begin{teo}[{\bf Hermitian Cauchy-Pompeiu formula in superspace}]\label{HC-PSSMat}
Let ${\bm G}_2^1$ be a matrix function of the form (\ref{CircMatGenForm}) with entries in $C^\infty(\Om)\otimes \mathfrak{G}_{2n}\otimes \mathcal{C}_{2m,2n}$, and let $g=g_0+{\bf g}\in C^\infty(\R^{2m})\otimes \mathfrak{G}_{2n}^{(ev)}$ be a phase function such that $\{g_0\leq 0\}\inc \Om$ is compact and $\pa_{\underline x}[g_0]\neq0$ on $g_0^{-1}(0)$. It then holds that
\begin{multline}\label{HC-PFormHerSupAnaMat}
\int_{\R^{{2m|2n}}} \hspace{-.4cm} {\bm \Psi}^{m|n}_{2\times 2} ({\bf Z}-{\bf U}) \,\del(g({\bf x}))\,{\bm{\mathcal D}}_{({\bf Z}, {\bf Z}^\dag)}[g({\bf x})] \, {\bm G}_2^1({\bf x}) -   \int_{\R^{2m|2n}}  \hspace{-.4cm} H(-g({\bf x})) \, {\bm \Psi}^{m|n}_{2\times 2} ({\bf Z}-{\bf U}) \left( {\bm{\mathcal D}}_{({\bf Z}, {\bf Z}^\dag)} {\bm G}_2^1({\bf x}) \right)\\= \begin{cases} {\bm G}_2^1({\bf y}), & g_0(\underline{y})<0, \\ 0,  &  g_0(\underline{y})>0.  \end{cases}
\end{multline}
\end{teo}
\pf
It is easily seen that ${\bm \Psi}^{m|n}_{2\times 2} ({\bf Z}-{\bf U})$ is a matrix distribution in the variable ${\bf Z}$ with singular support $\{\underline{y}\}$. Hence, for $g_0(\underline{y})\neq 0$ we have that 
\[sing \; supp \, H(-g{({\bf x})}) \cap sing \; supp \; {\bm \Psi}^{m|n}_{2\times 2} ({\bf Z}-{\bf U}) =\emptyset.\]
This means that {one can take} ${\bm {\Upsilon}}={\bm \Psi}^{m|n}_{2\times 2} ({\bf Z}-{\bf U})$ in (\ref{StokHeavmatSupSHerDist}). From Theorem \ref{FisrtMatApp} {it follows that}
\[{\bm \del}({\bf x}-{\bf y}) =  {\bm \Psi}^{m|n}_{2\times 2}({\bf Z}-{\bf U}) \; {\bm{\mathcal D}}_{({\bf Z}, {\bf Z}^\dag)},\]
{leading to}
\begin{multline}\label{StokHeavmatSupSHerDistInterm}
\int_{\R^{2m|2n}} \hspace{-.5cm} H(-g({\bf x})) \, {\del}({\bf x}-{\bf y}) {\bm G}_2^1({\bf x}) \\
   =  \int_{\R^{2m|2n}} \hspace{-.5cm} {\bm \Psi}^{m|n}_{2\times 2}({\bf Z}-{\bf U})  \,\del(g({\bf x}))  {\bm{\mathcal D}}_{({\bf Z}, {\bf Z}^\dag)}[g({\bf x})] \, {\bm G}_2^1({\bf x}) -  \int_{\R^{2m|2n}} \hspace{-.5cm} H(-g({\bf x})) \,  {\bm \Psi}^{m|n}_{2\times 2}({\bf Z}-{\bf U}) \, \left( {\bm{\mathcal D}}_{({\bf Z}, {\bf Z}^\dag)} {\bm G}_2^1({\bf x}) \right).
\end{multline}
Let us now examine the distributional product 
\[\del({\bf x}-{\bf y}) H(-g({\bf x}))= \frac{\pi^n}{n!}(\underline{x}\p-\underline{y}\p)^{2n}  \sum_{j=0}^n  \frac{(-{\bf g}({\bf x}))^j}{j!} \; \del(\underline{x}-\underline{y}) \del^{(j-1)}(-g_0({\underline x})).\]
It is clearly seen that $ sing \; supp \, \del(\underline{x}-\underline{y})=\{\underline{y}\}$ and  $ sing \; supp \, \del^{(j-1)}(-g_0({\underline x}))=g_0^{-1}(0)$, {whence}  (\ref{DistProd}) immediately shows for $g_0(\underline{y})\neq 0$ that 
\[\del(\underline{x}-\underline{y}) \del^{(j-1)}(-g_0({\underline x}))=0, \hspace{1cm} j=1,\ldots, n.\]
Thus $\del({\bf x}-{\bf y}) H(-g({\bf x}))=\del({\bf x}-{\bf y}) H(-g_0({\underline x}))$ and then, (\ref{DelProp}) yields
\begin{equation*}
\int_{\R^{2m|2n}} \del({\bf x}-{\bf y}) H(-g({\bf x})) {\bm G}_2^1({\bf x})=\int_{\R^{2m|2n}} \del({\bf x}-{\bf y}) H(-g_0({\underline x})) {\bm G}_2^1({\bf x})= H(-g_0({\underline y})) {\bm G}_2^1({\bf y}).
\end{equation*}
{Substitution of the later in (\ref{StokHeavmatSupSHerDistInterm}) gives the desired result (\ref{HC-PFormHerSupAnaMat}).}
$\hfill\square$

This theorem now leads to the following Hermitian Cauchy integral formulae in superspace. 
\begin{cor}
If the matrix function ${\bm G}_2^1$ is ${\bf SH}$-monogenic then,
\begin{equation*}
\int_{\R^{m|2n}} {\bm \Psi}^{m|n}_{2\times 2} ({\bf Z}-{\bf U}) \,\del(g({\bf x}))\,{\bm{\mathcal D}}_{({\bf Z}, {\bf Z}^\dag)}[g({\bf x})] \, {\bm G}_2^1({\bf x}) = \begin{cases} {\bm G}_2^1({\bf y}), & \hspace{-.1cm}g_0(\underline{y})<0, \\ 0,  &  \hspace{-.1cm} g_0(\underline{y})>0.  \end{cases}
\end{equation*}
\end{cor}

\begin{cor}\label{sh-monH-Cfor}
If the function $G$ is $sh$-monogenic then,
\begin{equation}\label{HCauchyFormHerSupAnaMatSpeCa}
\int_{\R^{m|2n}} {\bm \Psi}^{m|n}_{2\times 2} ({\bf Z}-{\bf U}) \,\del(g({\bf x}))\,{\bm{\mathcal D}}_{({\bf Z}, {\bf Z}^\dag)}[g({\bf x})] \, {\bm G}_0({\bf x}) = \begin{cases} {\bm G}_0({\bf y}), & \hspace{-.1cm}g_0(\underline{y})<0, \\ 0,  &  \hspace{-.1cm} g_0(\underline{y})>0.  \end{cases}
\end{equation}
\end{cor}

The above result may be considered as a Hermitian Cauchy integral theorem for the sh-monogenic function $G$. {For $n=0$ the above result} becomes the purely bosonic Hermitian Cauchy integral. The study of its boundary limits leads to Hermitian Clifford-Hardy spaces and to a Hermitian Hilbert transform, see e.g.\ \cite{MR2426334}.

\begin{remark}
The second summand {at} the left {hand} side of formula (\ref{HC-PFormHerSupAnaMat}) is the well-known {T\'eodorescu} transform, which is denoted by
\[
{\bm T}_g G_2^1({\bf y}) = -\int_{\R^{2m|2n}}  H(-g({\bf x})) \,  {\bm \Psi}^{m|n}_{2\times 2}({\bf Z}-{\bf U}) \,  {\bm G}_2^1({\bf x}).
\]
This operator constitutes a right inverse to the Dirac operator. Indeed, using Theorem \ref{FisrtMatApp} and (\ref{DelProp}) one easily obtains
\[
{\bm{\mathcal D}}_{({\bf U}, {\bf U}^\dag)} {\bm T}_g G_2^1({\bf y}) = \int_{\R^{2m|2n}} \hspace{-.5cm} H(-g({\bf x})) {\bf \del}({\bf x} -{\bf y}) {\bm G}_2^1({\bf x}) = \begin{cases} {\bm G}_2^1({\bf y}), & g_0(\underline{y})<0, \\ 0,  &  g_0(\underline{y})>0. \end{cases}
\]
The combination of the above inversion formula with the Hermitian Cauchy-Pompeiu Theorem \ref{HC-PSSMat}, yields a Hermitian Koppelman formula in superspace, see e.g.\ \cite{Brackx_SDh_MELE_Shap}.  This formula reads as follows for $g_0(\underline{y})<0$:
\begin{multline*}
\int_{\R^{m|2n}} \hspace{-.4cm} {\bm \Psi}^{m|n}_{2\times 2} ({\bf Z}-{\bf U}) \,\del(g({\bf x}))\,{\bm{\mathcal D}}_{({\bf Z}, {\bf Z}^\dag)}[g({\bf x})] \, {\bm G}_2^1({\bf x}) +  {\bm T}_g {\bm{\mathcal D}}_{({\bf Z}, {\bf Z}^\dag)} G_2^1({\bf y})
+{\bm{\mathcal D}}_{({\bf U}, {\bf U}^\dag)} {\bm T}_g G_2^1({\bf y})= 2{\bm G}_2^1({\bf y}).  
\end{multline*}
\end{remark}

\section{Integral formulae for holomorphic functions in superspace}\label{HolSSSection}

\subsection{Holomorphicity in superspace and sh-monogenicity}
Every superfunction $F({\bf x})\in C^\infty(\Om)\otimes \mathfrak{G}_{2n}$, $\Om\inc \R^{2m}$, can be written in terms of the Hermitian supervector variables ${\bf Z}$ and ${\bf Z}^\dag$ as
\begin{align}\label{SupFunFormCompSupVecVar}
F({\bf Z},{\bf Z}^\dag) &= \sum_{A,B\inc \{1,\ldots, n\}} F_{A,B} (\underline{Z}, \underline{Z}^\dag) \, \underline{Z}^{\p}_A \, {\underline{Z}^{\p}}^c_B, & F_{A,B}\in C^\infty(\Om),
\end{align}
where 
\begin{align*}
\underline{Z}^{\p}_A &= z\p_{j_1} \cdots z\p_{j_k}, & A&=\{j_1, \ldots, j_k\}, \quad 1\leq j_1< \ldots < j_k\leq n, 
\\
{\underline{Z}^{\p}}^c_B &=z\p_{\ell_1}^c \cdots z\p_{\ell_s}^c, & B&=\{\ell_1, \ldots, \ell_s\}, \quad 1\leq \ell_1< \ldots < \ell_s\leq n.
\end{align*}
The notion of {a} holomorphic function in the purely bosonic case {then} naturally extends to superfunctions.
\begin{defi}
A superfunction $F({\bf Z},{\bf Z}^\dag)$ of the form (\ref{SupFunFormCompSupVecVar}) is said to be holomorphic in the bosonic and fermionic complex variables $z_1,\ldots, z_m, z\p_1,\ldots, z\p_n$ if 
\begin{align*}
\pa_{z_j^c}[F] &= \pa_{z\p_k^c}[F]=0, & j&=1,\ldots, m, \quad k=1\ldots, n.
\end{align*}
\end{defi}
This holomorphicity  condition is equivalent to saying that the function $F$ does not depend on the conjugate variables $z_1^c,\ldots, z_m^c, z\p_1^c,\ldots, z\p_n^c$, i.e.\ 
\[F({\bf Z},{\bf Z}^\dag)= F({\bf Z}) = \sum_{A\inc \{1,\ldots, n\}} F_{A} (\underline{Z}) \, \underline{Z}^{\p}_A .\]

Let us now connect this  holomorphicity  notion in superspace with sh-monogenicity. To that end we start by introducing the classical primitive idempotent
\begin{align*}
I_b&=  \mathfrak f_1^\dag \mathfrak f_1 \cdots \mathfrak f_m^\dag \mathfrak f_m, &\mbox{ where } &&  \mathfrak f_j^\dag \mathfrak f_j&=\frac{1}{2}(1+ie_je_{m+j}), \quad j=1,\ldots, m.
\end{align*}
This idempotent clearly satisfies $e_j I_b=-ie_{m+j} I_b$ or equivalently  $\mathfrak f_j^\dag I_b=0$. 

A {similar} element to $I_b$ can be constructed in terms of the symplectic generators $e\p_j$'s. We first recall that the elements $e\p_{2j-1}$, $e\p_{2j}$ can be identified with the following operators when acting on the corresponding spinor space, see e.g.\ \cite{sommen2000extension, MR3375856},
\begin{align*}
e\p_{2j-1} &\fd e_{2m+1} \, \pa_{a_j}, & e\p_{2j} &\fd -e_{2m+1} \, a_j,
\end{align*}
where $a_1, \ldots, a_n$ are commuting variables and $e_{2m+1}$ is an {additional} orthogonal Clifford generator. It is readily seen that the above correspondences are in accordance with the commutation rules (\ref{CommRules}). We then need to find an element $I_f$ {similar} to $I_b$ in the spinor space that consists of all smooth functions in the variables $a_1, \ldots, a_n$. This element $I_f$ has to satisfy the key property 
\begin{align*}
e\p_{2j-1} I_f &=-i e\p_{2j} I_f, &\mbox{ or equivalently, } && \pa_{a_j} I_f=ia_j \,I_f.
\end{align*}
Such a function is given by 
\[I_f=\exp\, \left(\frac{i}{2}\sum_{j=1}^n a_j^2\right),\]
 which clearly is a null solution of the operator ${{\mathfrak f_j}\p}^\dag$, i.e.\
 \[{{\mathfrak f_j}\p}^\dag I_f=-\frac{1}{2}\left(e\p_{2j-1}+ie\p_{2j}\right)I_f=-\frac{e_{2m+1}}{2}\left(\pa_{a_j}-ia_j\right)[I_f]=0.\]
 \begin{pro}
A superfunction $F\in  C^\infty(\Om)\otimes \mathfrak{G}_{2n}$  is holomorphic in the variables $z_1,\ldots, z_m, z\p_1,\ldots, z\p_n$ if and only if the spinor-valued function $FI_b I_f$ is sh-monogenic.
\end{pro}
\pf
We first observe that 
\[\pa_{\bf Z} [FI_b I_f] = \sum_{j=1}^m \pa_{z_j}[F] \, \f_j^\dag \, I_b I_f + 2i\sum_{j=1}^n \pa_{z\p_j}[F]\, \f_j\p^\dag \, I_b I_f =0,
\]
{while} on the other hand,
\begin{equation}\label{PaZDagFI}
\pa_{{\bf Z}^\dag} [FI_b I_f] = \sum_{j=1}^m \pa_{z_j^c}[F] \, \f_j \, I_b I_f - 2i\sum_{j=1}^n \pa_{z\p_j^c}[F]\, \f_j\p \, I_b I_f ,
\end{equation}
{whence} it is clear that holomorphicity for $F$ implies sh-monogenicity for $FI_b I_f$.

\noindent Assume now that $FI_b I_f$ is sh-monogenic. In order to prove that $F$ is holomorphic, it suffices to {show} that all the elements $\f_j \, I_b I_f$, $\f_k\p \, I_b I_f$ are linearly independent when considering coefficients in $C^\infty(\Om) \otimes \mathfrak{G}_{2n}$, see (\ref{PaZDagFI}).

\noindent We {have},
\[
\f_j \, I_b =  \mathfrak f_1^\dag \mathfrak f_1  \cdots \left(\f_j \f_j^\dag\right) \f_j \cdots \mathfrak f_m^\dag \mathfrak f_m = \mathfrak f_1^\dag \mathfrak f_1  \cdots \left(1-\f_j^\dag \f_j \right) \f_j \cdots \mathfrak f_m^\dag \mathfrak f_m = \mathfrak f_1^\dag \mathfrak f_1  \cdots \mathfrak f_{j-1}^\dag \mathfrak f_{j-1} \,  \f_j \, \mathfrak f_{j+1}^\dag \mathfrak f_{j+1} \cdots \mathfrak f_m^\dag \mathfrak f_m,
\]
{yielding},
\begin{equation}\label{MultPrimIdemBos}
\f_k^\dag \f_j I_b=\del_{k,j} I_b.
\end{equation}
Moreover,
\[
\f_j\p I_b I_f = I_b \f_j\p I_f=I_b \frac{e_{2m+1}}{2} (\pa_{a_j}+i a_j) \, \exp \left(\frac{i}{2}\sum_{j=1}^n a_j^2\right)= I_b e_{2m+1} i a_j \exp \left(\frac{i}{2}\sum_{j=1}^n a_j^2\right).
\]
Hence, taking into account (\ref{PaZDagFI}), the sh-monogenicity of $FI_b I_f$ reduces to  
\[ \left(\sum_{j=1}^m \pa_{z_j^c}[F] \, \f_j \, I_b - 2i\sum_{j=1}^n \pa_{z\p_j^c}[F]\, I_b e_{2m+1} i a_j\right) \exp \left(\frac{i}{2}\sum_{j=1}^n a_j^2\right) =0,\]
{implying}
\[ \sum_{j=1}^m \pa_{z_j^c}[F] \, \f_j \, I_b - 2i\sum_{j=1}^n \pa_{z\p_j^c}[F]\, I_b e_{2m+1} i a_j =0.\]
Multiplying {the above expression} from the left by $\f_k^\dag$ we get, on account of (\ref{MultPrimIdemBos}), that $\pa_{z_k^c}[F] \, I_b=0$. This directly implies that $\pa_{z_k^c}[F]=0$ for every $k=1,\ldots, m$, {whence the previous equality reduces to} 
\[\sum_{j=1}^n \pa_{z\p_j^c}[F]\, I_b e_{2m+1} i a_j =0.\]
Finally, by taking $a_k=1$ and $a_j=0$ ($j\neq k$) we obtain $\pa_{z\p_k^c}[F]=0$ for every  $k=1,\ldots, n$.
$\hfill\square$

\subsection{Bochner-Martinelli theorem for holomorphic superfunctions}
The above result shows that considering functions of the form $FI_bI_f$ establishes a connection between Hermitian {monogenicity} and holomorphic functions in superspace. In this section, we investigate the nature of the Hermitian Cauchy integral formula obtained in Corollary \ref{sh-monH-Cfor} for this type of functions.

To this end, we will explicitly compute the left-hand side of (\ref{HCauchyFormHerSupAnaMatSpeCa}), taking $G=FI_bI_f$ where $F({\bf Z}) =\sum_{A} F_{A} (\underline{Z}) \, \underline{Z}^{\p}_A$ is a holomorphic function. We first obtain,
\[
{\bm{\mathcal D}}_{({\bf Z}, {\bf Z}^\dag)}[g({\bf x})] \, {\bm G}_0({\bf x}) = \left(\begin{array}{cc} \pa_{{\bf Z}}[g({\bf x})] & * \\ \pa_{{\bf Z}^\dag}[g({\bf x})] & * \end{array}\right) \left(\begin{array}{cc} F({\bf Z})I_bI_f & * \\ 0 & * \end{array}\right) =
\left(\begin{array}{cc} \pa_{{\bf Z}}[g({\bf x})] F({\bf Z})I_bI_f & * \\ \pa_{{\bf Z}^\dag}[g({\bf x})] F({\bf Z})I_bI_f & * \end{array}\right),
\]
where the second columns have not been written, since {they} only duplicate the first ones (in reversed order) on account of the circulant structure of the involved matrices. Further calculation yields,
\[
\pa_{{\bf Z}}[g({\bf x})] F({\bf Z})I_bI_f =  \sum_{j=1}^m \pa_{z_j}[g({\bf x})] \, F({\bf Z}) \, (\f_j^\dag \, I_b) I_f + 2i\sum_{j=1}^n \pa_{z\p_j}[g({\bf x})]\, F({\bf Z}) \, I_b \, (\f_j\p^\dag \, I_f)=0,
\]
and as a consequence, ${\bm{\mathcal D}}_{({\bf Z}, {\bf Z}^\dag)}[g({\bf x})] \, {\bm G}_0({\bf x}) = \left(\begin{array}{cc} 0 & * \\ \pa_{{\bf Z}^\dag}[g({\bf x})] F({\bf Z})I_bI_f & * \end{array}\right)$. 
Hence
\begin{align*}
{\bm \Psi}^{m|n}_{2\times 2} ({\bf Z}-{\bf U}) \,\del(g({\bf x}))\,{\bm{\mathcal D}}_{({\bf Z}, {\bf Z}^\dag)}[g({\bf x})] \, {\bm G}_0({\bf x})  &= \left(\begin{array}{cc}{\Psi^{m|n}_1} ({\bf Z}-{\bf U}) &* \\ {\Psi^{m|n}_1}^\dag ({\bf Z}-{\bf U}) & * \end{array}\right)  \left(\begin{array}{cc} 0 & * \\ \del(g({\bf x}))\, \pa_{{\bf Z}^\dag}[g({\bf x})] \; F({\bf Z})I_bI_f & * \end{array}\right) \\ 
\\
&=  \left(\begin{array}{cc}  {\Psi^{m|n}_1}^\dag ({\bf Z}-{\bf U})\del(g({\bf x}))\, \pa_{{\bf Z}^\dag}[g({\bf x})] \; F({\bf Z})I_bI_f  & * \\ {\Psi^{m|n}_1} ({\bf Z}-{\bf U})  \,\del(g({\bf x}))\, \pa_{{\bf Z}^\dag}[g({\bf x})] \; F({\bf Z})I_bI_f & * \end{array}\right).
\end{align*}
Denoting the even element $|{\bf Z}-{\bf U}|={\bm \rho}$ we {compute}
\begin{multline*}
{\Psi^{m|n}_1}^\dag ({\bf Z}-{\bf U})\del(g({\bf x}))\, \pa_{{\bf Z}^\dag}[g({\bf x})]
 \\= \frac{2\del(g({\bf x}))}{|\Sa^{2m-1|2n}| {\bm \rho}^M }  \hspace{-.1cm} \left[ \left(\underline{Z}- \underline{U}\right)^\dag \pa_{{\underline{Z}}^\dag}[g({\bf x})] \hspace{-.05cm} +  \hspace{-.05cm} \left(\underline{Z}- \underline{U}\right)^\dag \pa_{{\underline{Z}^{\p}}^\dag}[g({\bf x})] \hspace{-.05cm} +\hspace{-.05cm} \left(\underline{Z}^{\p}- \underline{U}^{\p}\right)^\dag \hspace{-.05cm} \pa_{{\underline{Z}^{\p}}^\dag}[g({\bf x})] \hspace{-.05cm} + \hspace{-.05cm} \left(\underline{Z}^{\p}- \underline{U}^{\p}\right)^\dag \hspace{-.05cm} \pa_{{\underline{Z}}^\dag}[g({\bf x})] \right].  
\end{multline*}
{We now consider each term in the previous sum separately, obtaining}
\begin{align*}
\left(\underline{Z}- \underline{U}\right)^\dag \pa_{{\underline{Z}}^\dag}[g({\bf x})] &= \left(\sum_{j=1}^m (z_j-u_j)^c \, \f_j^\dag\right)   \left(\sum_{k=1}^m \f_k \, \pa_{z_k^c}[g({\bf x})]\right) \\
&= - \sum_{j\neq k}  (z_j-u_j)^c \, \pa_{z_k^c}[g({\bf x})] \, \f_k \f_j^\dag + \sum_{j=1}^m (z_j-u_j)^c \, \pa_{z_j^c}[g({\bf x})] \, (1-\f_j \f_j^\dag ), \\
\left(\underline{Z}- \underline{U}\right)^\dag \pa_{{\underline{Z}^{\,\p}}^\dag}[g({\bf x})] &= -2i \left(\sum_{j=1}^m (z_j-u_j)^c \, \f_j^\dag\right)   \left(\sum_{k=1}^n \f_k\p \, \pa_{z\p_k^c}[g({\bf x})]\right) = 2i \hspace{-.2cm}\sum_{1\le j\le m\atop 1\le k\le n} (z_j-u_j)^c \, \pa_{z\p_k^c}[g({\bf x})] \, \f_k\p \, \f_j^\dag, \\
\left(\underline{Z}^{\p}- \underline{U}^{\p}\right)^\dag \pa_{{\underline{Z}^{\,\p}}^\dag}[g({\bf x})] &= -2i\left(\sum_{j=1}^n (z\p_j-u\p_j)^c \, \f_j\p^\dag\right)   \left(\sum_{k=1}^n \f_k\p \, \pa_{z\p_k^c}[g({\bf x})]\right) \\
&= -2i \sum_{j\neq k}  (z\p_j-u\p_j)^c \, \pa_{z\p_k^c}[g({\bf x})] \, \f_k\p \, \f_j\p^\dag     -2i \sum_{j=1}^n  (z\p_j-u\p_j)^c \, \pa_{z\p_j^c}[g({\bf x})] \,\left(\frac{i}{2}+\f_j\p \, {\f_j\p}^\dag \right),\\
\left(\underline{Z}^{\p}- \underline{U}^{\p}\right)^\dag \pa_{{\underline{Z}}^\dag}[g({\bf x})] &= 
 \left(\sum_{j=1}^n (z\p_j-u\p_j)^c \, \f_j\p^\dag\right)   \left(\sum_{k=1}^m \f_k \, \pa_{z_k^c}[g({\bf x})]\right) = -\hspace{-.2cm} \sum_{1\le j\le n\atop 1\le k\le m} (z\p_j-u\p_j)^c \, \pa_{z_k^c}[g({\bf x})] \, \f_k \f_j\p^\dag.
 \end{align*}
This yields,
 \begin{align*}
 \left(\underline{Z}- \underline{U}\right)^\dag \pa_{{\underline{Z}}^\dag}[g({\bf x})] \,  F({\bf Z})I_bI_f &= \sum_{j=1}^m (z_j-u_j)^c \, \pa_{z_j^c}[g({\bf x})] \, F({\bf Z})I_bI_f, & 
 \left(\underline{Z}- \underline{U}\right)^\dag \pa_{{\underline{Z}^{\p}}^\dag}[g({\bf x})]\,  F({\bf Z})I_bI_f &= 0, \\
  \left(\underline{Z}^{\p}- \underline{U}^{\p}\right)^\dag \pa_{{\underline{Z}^{\p}}^\dag}[g({\bf x})]  \,  F({\bf Z})I_bI_f  &= \sum_{j=1}^n  (z\p_j-u\p_j)^c \, \pa_{z\p_j^c}[g({\bf x})]  \,  F({\bf Z})I_bI_f,
& 
 \left(\underline{Z}^{\p}- \underline{U}^{\p}\right)^\dag \pa_{{\underline{Z}}^\dag}[g({\bf x})]  \,  F({\bf Z})I_bI_f  &= 0 .
 \end{align*}
 Hence we obtain 
  \begin{align*}
   &{\Psi^{m|n}_1}^\dag ({\bf Z}-{\bf U})\del(g({\bf x}))\, \pa_{{\bf Z}^\dag}[g({\bf x})]\,  F({\bf Z})I_bI_f   \\
   &\quad \quad = \frac{2\del(g({\bf x}))}{|\Sa^{2m-1|2n}| {\bm \rho}^M } \left[ \sum_{j=1}^m (z_j-u_j)^c \, \pa_{z_j^c}[g({\bf x})] +
\sum_{j=1}^n  (z\p_j-u\p_j)^c \, \pa_{z\p_j^c}[g({\bf x})] \right] \,  F({\bf Z})I_bI_f  \\
& \quad \quad = \frac{2\del(g({\bf x}))}{|\Sa^{2m-1|2n}| {\bm \rho}^M } \, D^\dag_{{\bf Z}-{\bf U}, {\bf Z}} [g({\bf x})] \,   F({\bf Z})I_bI_f, 
  \end{align*}
 where the differential operator
 \[
 D^\dag_{{\bf Z}-{\bf U}, {\bf Z}}=  \sum_{j=1}^m (z_j-u_j)^c \, \pa_{z_j^c} +
\sum_{j=1}^n  (z\p_j-u\p_j)^c \, \pa_{z\p_j^c} = \left\{ \pa_{\bf Z}^\dag, ({\bf Z}-{\bf U})^\dag\right\} -\frac{1}{2}\left(\frac{M}{2}-i{\bf B}\right)
\]
is the so-called directional derivative with respect to ${\bf Z}^\dag$ in the direction $({\bf Z}-{\bf U})^\dag$, see \cite{DS_Guz_Somm2}.

Thus, the Hermitian Cauchy integral formula (\ref{HCauchyFormHerSupAnaMatSpeCa}) for $g_0(\underline{y})<0$ yields the following two statements:
\begin{align}
 \frac{2}{|\Sa^{2m-1|2n}| } \int_{\R^{2m|2n}} \frac{\del(g({\bf x}))}{|{\bf Z}-{\bf U}|^M} \, D^\dag_{{\bf Z}-{\bf U}, {\bf Z}} [g({\bf x})] \,   F({\bf Z})I_bI_f  &= F({\bf U})I_bI_f ,\nonumber\\
  \frac{2}{|\Sa^{2m-1|2n}| } \int_{\R^{2m|2n}} \frac{\del(g({\bf x}))}{|{\bf Z}-{\bf U}|^M} \,({\bf Z}-{\bf U}) \, \pa_{{\bf Z}^\dag}[g({\bf x})]\,  F({\bf Z})I_bI_f&= 0. \label{B-Mcons2}
\end{align}
The first {one}  leads to the following integral representation of holomorphic superfunctions. 
\begin{teo}[{\bf Bochner-Martinelli formula in superspace}]\label{B-MDistSSTheom}
Let $F({\bf Z})\in C^\infty(\Om)\otimes \mathfrak{G}_{2n}$ be a holomorphic function in the variables $z_1\ldots, z_m, z\p_1, \ldots, z\p_n$, ($m>n$), and let $g=g_0+{\bf g}\in C^\infty(\R^{2m})\otimes \mathfrak{G}_{2n}^{(ev)}$ be a phase function such that $\{g_0\leq 0\}\inc \Om$ is compact and $\pa_{\underline x}[g_0]\neq0$ on $g_0^{-1}(0)$. It then follows for $g_0(\underline y)<0$ that 
\begin{equation}\label{SupBoc-MartForm}
 \frac{2}{|\Sa^{2m-1|2n}| } \int_{\R^{2m|2n}} \frac{\del(g({\bf x}))}{|{\bf Z}-{\bf U}|^M} \, D^\dag_{{\bf Z}-{\bf U}, {\bf Z}} [g({\bf x})] \,   F({\bf Z}) = F({\bf U}).
\end{equation}
\end{teo}
\begin{remark}
The above theorem {indeed constitutes} an extension of the classical Bochner-Martinelli formula to superspace. In the next section it will be shown that (\ref{SupBoc-MartForm}) reduces to (\ref{B-Moriginal}) when $n=0$.
\end{remark}
On the other hand, the second statement yields the following result. 
\begin{teo}
Under the same conditions {as in} Theorem  \ref{B-MDistSSTheom} one has
\begin{align*}
\int_{\R^{2m|2n}} \hspace{-.2cm}  \del(g({\bf x}))\, \frac{ z_j-u_j}{|{\bf Z}-{\bf U}|^M} \,\pa_{z_k^c}[g({\bf x})] \, F({\bf Z}) &= \phantom{-} \int_{\R^{2m|2n}} \hspace{-.2cm} \del(g({\bf x}))\, \frac{ z_k-u_k}{|{\bf Z}-{\bf U}|^M} \,\pa_{z_j^c}[g({\bf x})] \, F({\bf Z}), \\
-2i \int_{\R^{2m|2n}} \hspace{-.2cm} \del(g({\bf x}))\, \frac{ z_j-u_j}{|{\bf Z}-{\bf U}|^M} \,\pa_{z\p_k^c}[g({\bf x})] \, F({\bf Z}) &= \phantom{-} \int_{\R^{2m|2n}} \hspace{-.2cm} \del(g({\bf x}))\, \frac{ z\p_k-u\p_k}{|{\bf Z}-{\bf U}|^M} \,\pa_{z_j^c}[g({\bf x})] \, F({\bf Z}), \\
\phantom{-2i}  \int_{\R^{2m|2n}} \hspace{-.2cm} \del(g({\bf x}))\, \frac{ z\p_j-u\p_j}{|{\bf Z}-{\bf U}|^M} \,\pa_{z\p_k^c}[g({\bf x})] \, F({\bf Z}) &= -\int_{\R^{2m|2n}} \hspace{-.2cm} \del(g({\bf x}))\, \frac{ z\p_k-u\p_k}{|{\bf Z}-{\bf U}|^M} \,\pa_{z\p_j^c}[g({\bf x})] \, F({\bf Z}). 
\end{align*}
\end{teo}
\pf
The proof directly follows from expanding expression (\ref{B-Mcons2}). We {have}
\[
({\bf Z}-{\bf U}) \, \pa_{{\bf Z}^\dag}[g({\bf x})] = \left(\underline{Z}- \underline{U}\right) \pa_{{\underline{Z}}^\dag}[g({\bf x})] +
\left(\underline{Z}- \underline{U}\right) \pa_{{\underline{Z}^{\p}}^\dag}[g({\bf x})]+\left(\underline{Z}^{\p}- \underline{U}^{\p}\right) \pa_{{\underline{Z}}^\dag}[g({\bf x})]+\left(\underline{Z}^{\p}- \underline{U}^{\p}\right) \pa_{{\underline{Z}^{\p}}^\dag}[g({\bf x})],
\]
where
\begin{align*}
\left(\underline{Z}- \underline{U}\right) \pa_{{\underline{Z}}^\dag}[g({\bf x})] &= \phantom{-2i}  \left(\sum_{j=1}^m (z_j-u_j) \, \f_j\right)   \left(\sum_{k=1}^m \f_k \, \pa_{z_k^c}[g({\bf x})]\right) = \phantom{-2i} \sum_{j\neq k}  (z_j-u_j) \, \pa_{z_k^c}[g({\bf x})] \, \f_j \f_k, \\
\left(\underline{Z}- \underline{U}\right) \pa_{{\underline{Z}^{\p}}^\dag}[g({\bf x})] &= -2i \left(\sum_{j=1}^m (z_j-u_j) \, \f_j\right)   \left(\sum_{k=1}^n \f_k\p \, \pa_{z\p_k^c}[g({\bf x})]\right) = -2i \sum_{1\le j\le m\atop 1\le k\le n} (z_j-u_j) \, \pa_{z\p_k^c}[g({\bf x})] \, \f_j \f_k\p , \\
\left(\underline{Z}^{\p}- \underline{U}^{\p}\right) \pa_{{\underline{Z}}^\dag}[g({\bf x})] &= 
 \phantom{-2i}  \left(\sum_{k=1}^n (z\p_k-u\p_k) \, \f_k\p\right)   \left(\sum_{j=1}^m \f_j \, \pa_{z_j^c}[g({\bf x})]\right) = \phantom{-2i} \sum_{1\le j\le m\atop 1\le k\le n} (z\p_k-u\p_k) \, \pa_{z_j^c}[g({\bf x})] \, \f_k\p \, \f_j, \\
\left(\underline{Z}^{\p}- \underline{U}^{\p}\right) \pa_{{\underline{Z}^{\p}}^\dag}[g({\bf x})] &= -2i\left(\sum_{j=1}^n (z\p_j-u\p_j) \, \f_j\p\right)   \left(\sum_{k=1}^n \f_k\p \, \pa_{z\p_k^c}[g({\bf x})]\right) 
= -2i \sum_{1\leq j,k\leq n}  (z\p_j-u\p_j) \, \pa_{z\p_k^c}[g({\bf x})] \, \f_j\p \, \f_k\p.
\end{align*}
Hence (\ref{B-Mcons2}) reads
\begin{multline*}
 \frac{2}{|\Sa^{2m-1|2n}| }  \int_{\R^{2m|2n}} \frac{\del(g({\bf x}))}{|{\bf Z}-{\bf U}|^M}  \sum_{1\leq j< k\leq m}  \left[ (z_j-u_j) \, \pa_{z_k^c}[g({\bf x})] - (z_k-u_k) \, \pa_{z_j^c}[g({\bf x})] \right] F({\bf Z}) \, \f_j \f_k \, I_b I_f \\
 + \frac{2}{|\Sa^{2m-1|2n}| }  \int_{\R^{2m|2n}} \frac{\del(g({\bf x}))}{|{\bf Z}-{\bf U}|^M}  \sum_{1\le j\le m\atop 1\le k\le n}  \left[ -2i(z_j-u_j) \, \pa_{z\p_k^c}[g({\bf x})] - (z\p_k-u\p_k) \, \pa_{z_j^c}[g({\bf x})] \right] F({\bf Z}) \, \f_j \f_k\p \, I_b I_f \\
 - \frac{4i}{|\Sa^{2m-1|2n}| }  \int_{\R^{2m|2n}} \frac{\del(g({\bf x}))}{|{\bf Z}-{\bf U}|^M}  \sum_{1\leq j< k\leq n}  \left[ (z\p_j-u\p_j) \, \pa_{z\p_k^c}[g({\bf x})] + (z\p_k-u\p_k) \, \pa_{z\p_j^c}[g({\bf x})] \right] F({\bf Z}) \, \f_j\p \, \f_k\p \, I_b I_f \\
 - \frac{4i}{|\Sa^{2m-1|2n}| }  \int_{\R^{2m|2n}} \frac{\del(g({\bf x}))}{|{\bf Z}-{\bf U}|^M}  \sum_{1\leq j \leq n}  (z\p_j-u\p_j) \, \pa_{z\p_j^c}[g({\bf x})]\,F({\bf Z}) \, {\f_j\p}^{\,2}  \, I_b I_f =0.
\end{multline*}
Thus, it suffices to prove {that} all the elements 
\begin{align}\label{WeirElem}
 \f_j \f_k \, I_b I_f, & \quad   1\leq j< k\leq m, &  \f_j \f_k\p \, I_b I_f,  &\quad 1\le j\le m, \;\; 1\le k\le n,  & \f_j\p \, \f_k\p \, I_b I_f &\quad 1\leq j\leq k\leq n,
 \end{align}
 are linearly independent when considering coefficients in $C^\infty(\Om)\otimes \mathfrak{G}_{2n}$. {So take} $A_{j,k}$, $B_{j,k}$, $C_{j,k}$ $\in C^\infty(\Om)\otimes \mathfrak{G}_{2n}$ such that
 \begin{equation}\label{LIWeirdElem}
 \sum_{1\leq j< k\leq m} A_{j,k}\, \f_j \f_k \, I_b I_f+  \sum_{1\le j\le m\atop 1\le k\le n} B_{j,k} \, \f_j \f_k\p \, I_b I_f+ \sum_{1\leq j\leq k\leq n}  C_{j,k} \, \f_j\p \, \f_k\p \, I_b I_f =0.
 \end{equation}
 We now observe that
 \begin{align*}
 \f_\ell^\dag \,  \f_s^\dag \left(  \f_j \f_k \, I_b I_f\right) &= (\del_{s,j} \del_{\ell,k} - \del_{\ell,j}\del_{s,k}) I_bI_f, &  \f_\ell^\dag \,  \f_s^\dag \left( \f_j \f_k\p \, I_b I_f\right) &=0,  &  \f_\ell^\dag \,  \f_s^\dag \left(  \f_j\p \, \f_k\p \, I_b I_f\right) &=0.
 \end{align*}
{Multiplying} (\ref{LIWeirdElem}) from the left by $ \f_\ell^\dag \,  \f_s^\dag$ ($1\leq s < \ell \leq m$) we obtain $A_{s,\ell}\, I_b I_f=0$ which implies  $A_{s,\ell}=0$, {whence we are left with} 
\[
\sum_{1\le j\le m\atop 1\le k\le n} B_{j,k} \, \f_j \f_k\p \, I_b I_f+ \sum_{1\leq j\leq k\leq n}  C_{j,k} \, \f_j\p \, \f_k\p \, I_b I_f =0.
\]
In the same order of ideas we get,
\begin{align*}
\f_\ell\p^\dag \, \f_s^\dag \, ( \f_j \f_k\p \, I_b I_f) &= \frac{i}{2} \del_{\ell,k} \del_{s,j} I_bI_f, & \f_\ell\p^\dag \, \f_s^\dag \, ( \f_j\p \, \f_k\p \, I_b I_f ) &=0.
\end{align*}
{whence multiplying the remainder of} (\ref{LIWeirdElem}) from the left by $\f_\ell\p^\dag \, \f_s^\dag$ ($1\le s\le m$,  $1\le \ell\le n$) {yields} $B_{s,\ell}=0$, {thus further reducing the equality to} 
\[
\sum_{1\leq j\leq k\leq n}  C_{j,k} \, \f_j\p \, \f_k\p \, I_b I_f =0.
\]
Finally, we compute
\begin{align*}
\f_\ell\p^\dag \, \f_s\p^\dag \, ( \f_j\p \, \f_k\p \, I_b I_f) = -\frac{1}{4} (\del_{\ell,j} \del_{s,k} + \del_{s,j}\del_{\ell,k}) I_bI_f,
\end{align*}
which allows to {conclude that} $C_{\ell,s}=0$ after multiplying (\ref{LIWeirdElem}) by  $\f_\ell\p^\dag \, \f_s\p^\dag$ ($1\leq \ell\leq s\leq n$). In this way, we have proven that all {coefficients in (\ref{LIWeirdElem}) are zero, meaning that all} elements (\ref{WeirElem}) {indeed} are linearly independent. 
$\hfill\square$

\subsection{Some examples}
In this section we study some particular but important applications of the Bochner-Martinelli formula in superspace. 
\\[+.2cm]
\noindent {\bf Case 1.} 

\noindent We first consider the case of a purely bosonic phase function $g({\bf x})=g_0(\underline{x}) \in C^\infty(\R^{2m})$, i.e.\ ${\bf g}=0$, which satisfies the same conditions as in Theorem \ref{B-MDistSSTheom}. In this case, formula (\ref{SupBoc-MartForm}) reads 
\begin{align}\label{g=g_0B-M}
 \frac{2}{|\Sa^{2m-1|2n}| } \int_{\R^{2m|2n}} \frac{\del(g_0(\underline{x}))}{|{\bf Z}-{\bf U}|^M} \, \left( \sum_{j=1}^m (z_j-u_j)^c \pa_{z_j^c}[g_0({\underline x})]\right) \,   F({\bf Z}) &= F({\bf U}),  & g_0(\underline y)&<0.
\end{align}
This formula reduces to the classical Bochner-Martinelli formula (\ref{B-Moriginal}) as we will show next. We begin by recalling the following classical result for surface integration over $\Gam:=g_0^{-1}(0)$, see \cite[p.~136]{MR1996773},
\[
\int_{\R^{2m}} \del(g_0(\underline{x})) |\pa_{\underline{x}}[g_0(\underline{x})]| f(\underline{x}) \, dV_{\underline{x}} =\int_{\Gam} f(\underline{x})\, dS_{\underline{x}}.
\]
The $j$-th coordinate $n_j(\underline{x})$ of the exterior normal vector $n(\underline{x})$ to the surface $\Gam$ at the point $\underline{x}\in \Gam$ is given by $n_j(\underline{x})= \frac{\pa_{x_j}[g](\underline{x})}{ |\pa_{\underline{x}}[g_0(\underline{x})]|}$. Hence, from the above formula one easily obtains
\[
\int_{\R^{2m}} \del(g_0(\underline{x})) \pa_{x_j}[g_0(\underline{x})] f(\underline{x}) \, dV_{\underline{x}} = \int_{\Gam} n_j(\underline{x})f(\underline{x})\, dS_{\underline{x}}.
\]
Moreover, since $\Gam$ is a smooth surface in $\R^{2m}$, we can write $\widehat{d x_j}= (-1)^{j-1} n_j(\underline{x}) \, dS_{\underline{x}}$ where $\widehat{d x_j}$ is the differential ($2m-1$)-form
 \[
 \widehat{d x_j} = dx_1 \wedge\cdots \wedge dx_{j-1} \wedge d x_{j+1} \wedge\ldots\wedge dx_{2m}.
 \]
This allows to change the above distributional approach to classical surface integration by differential forms. In particular, we have that
\begin{align*} 
\int_{\R^{2m}} \del{(g_0(\underline{x}))} \, \pa_{z_j^c}[g_0(\underline{x})] f(\underline{x}) \, dV_{\underline{x}} &=
 \frac{1}{2} \int_{\R^{2m}} \del{(g_0(\underline{x}))} \left(\pa_{x_j}+i \pa_{x_{m+j}}\right)\hspace{-.1cm}[g_0(\underline{x})] \, f(\underline{x}) \, dV_{\underline{x}} \\
 &= \frac{1}{2} \int_{\Gam} \big( n_j(\underline{x}) + i n_{m+j}(\underline{x})\big) f(\underline{x}) \, dS_{\underline{x}} \\
 &=  \frac{1}{2} \int_{\Gam}  \left( (-1)^{j-1}  \widehat{d x_j} +i (-1)^{m+j-1}  \widehat{d x_{m+j}} \right)  f(\underline{x})  .
\end{align*}
We now write,
\begin{align*}
 (-1)^{j-1}  \widehat{d x_j}  &= (-1)^{\frac{m(m-1)}{2}} \onda{\widehat{d x_j}}, &  -(-1)^{m+j-1}  \widehat{d x_{m+j}}  &= (-1)^{\frac{m(m-1)}{2}} \onda{\widehat{d x_{m+j}}},
\end{align*}
with
\begin{align*}
\onda{\widehat{d x_j}} &= \left(dx_1\wedge dx_{m+1}\right)\wedge \ldots \wedge \left([dx_{j}]\wedge dx_{m+j}\right) \wedge \ldots \wedge \left(dx_m\wedge dx_{2m}\right), & j&=1,\ldots, m, \\
\onda{\widehat{d x_{m+j}}} &= \left(dx_1\wedge dx_{m+1}\right)\wedge \ldots \wedge \left(dx_{j}\wedge [dx_{m+j}]\right) \wedge \ldots \wedge \left(dx_m\wedge dx_{2m}\right), & j&=1,\ldots, m,
\end{align*} 
where $[\cdot]$ denotes omitting that particular differential. We then obtain, 
\[
\int_{\R^{2m}} \del{(g_0(\underline{x}))} \, \pa_{z_j^c}[g_0(\underline{x})] f(\underline{x}) \, dV_{\underline{x}} = \frac{(-1)^{\frac{m(m-1)}{2}}}{2} \int_{\Gam}  \left(   \onda{\widehat{d x_j}} -i \onda{\widehat{d x_{m+j}}} \right)  f(\underline{x}). \, 
\]
On the other hand, 
\[ \onda{\widehat{d x_j}} -i \onda{\widehat{d x_{m+j}}} = -2 \left(\frac{i}{2}\right)^m \widehat{d z_j^c} =  (-1)^{\frac{m(m+1)}{2}-j} (-2) \left(\frac{i}{2}\right)^m  \onda{\widehat{d z_j^c}}\]
where we have introduced the complex differential forms
\begin{align*}
\widehat{d z_j^c} &=  \left(dz_1\wedge dz^c_{1}\right)\wedge\ldots \wedge \left(dz_{j}\wedge [dz^c_{j}]\right) \wedge \ldots \wedge \left(dz_m\wedge dz^c_{m}\right), \\
\onda{\widehat{d z_j^c}} &= dz^c_{1} \wedge \ldots \wedge  [dz^c_{j}]  \wedge \ldots \wedge dz^c_{m} \wedge dz_1\wedge \ldots \wedge dz_m,
\end{align*}
which are clearly connected by $\widehat{d z_j^c} = (-1)^{\frac{m(m+1)}{2}-j} \;  \onda{\widehat{d z_j^c}}$.

Hence, 
\begin{align*}
\int_{\R^{2m}} \del{(g_0(\underline{x}))} \, \pa_{z_j^c}[g_0(\underline{x})] f(\underline{x}) \, dV_{\underline{x}} = (-1)^{m-j+1} \left(\frac{i}{2}\right)^m \int_\Gam  f(\underline{x}) \, \onda{\widehat{d z_j^c}}.
\end{align*}
Applying this last result to the left side of formula (\ref{g=g_0B-M}) yields
 \begin{align*}
 & \frac{2}{|\Sa^{2m-1|2n}| } \int_B \sum_{j=1}^m \left( \int_{\R^{2m}}  \del{(g_0(\underline{x}))} \, \pa_{z_j^c}[g_0(\underline{x})]  \frac{(z_j-u_j)^c}{|{\bf Z}-{\bf U}|^M}  \,   F({\bf Z}) \, dV_{\underline{x}} \right) \\
& \quad \quad\quad=\frac{2}{|\Sa^{2m-1|2n}| } \int_B \sum_{j=1}^m  (-1)^{m-j+1} \left(\frac{i}{2}\right)^m  \int_\Gam \, \onda{\widehat{d z_j^c}} \frac{(z_j-u_j)^c}{|{\bf Z}-{\bf U}|^M}  \,   F({\bf Z}) \\
&\quad \quad\quad =\frac{2}{|\Sa^{2m-1|2n}| (2i)^m}  \int_\Gam  \int_B \sum_{j=1}^m (-1)^{j-1}  \onda{\widehat{d z_j^c}} \frac{(z_j-u_j)^c}{|{\bf Z}-{\bf U}|^M}  \,   F({\bf Z}).
\end{align*}
Then  (\ref{g=g_0B-M}) can be rewritten as
\begin{align*}
\frac{(m-n-1)! \, \pi^n}{ (2i\pi)^m}  \int_\Gam  \int_B \left(  \sum_{j=1}^m (-1)^{j-1}  \onda{\widehat{d z_j^c}} \frac{(z_j-u_j)^c}{|{\bf Z}-{\bf U}|^M}\right) \, F({\bf Z}) &= F({\bf U}),  & g_0(\underline y)&<0,
\end{align*}
which exactly coincides with formula (\ref{B-Moriginal}) for $n=0$.
\\[+.2cm]

\noindent{\bf Case 2.}

 \noindent We now examine the form {which} (\ref{SupBoc-MartForm}) takes on the supersphere of radius $R>0$ defined by means of the phase function $g({\bf x}) = -{\bf x}^2 -R^2$. Observe that
\[g({\bf x}) = |{\bf x}|^2 -R^2 = |{\bf Z}|^2 -R^2=\{{\bf Z}, {\bf Z}^\dag\}  -R^2=\sum_{j=1}^m z_jz_j^c-\frac{i}{2}\sum_{j=1}^n z\p_jz\p_j^c -R^2.\]
Then, $\pa_{z_j^c}[g({\bf x})]= z_j$ and  $\pa_{z\p_j^c}[g({\bf x})]= \frac{i}{2}z\p_j$, {leading to}
\[
 D^\dag_{{\bf Z}-{\bf U}, {\bf Z}} [g({\bf x})] =  \sum_{j=1}^m (z_j-u_j)^c \, z_j +\frac{i}{2} \sum_{j=1}^n  (z\p_j-u\p_j)^c \, z\p_j = \{{\bf Z}, ({\bf Z}-{\bf U})^\dag \}.
\]
Hence, the Bochner-Martinelli formula on the supersphere of radius $R>0$ takes the form 
\begin{align*}
\frac{2}{|\Sa^{2m-1|2n}|}  \int_{\R^{2m|2n}} \frac{\del(\{{\bf Z}, {\bf Z}^\dag\}  -R^2)}{|{\bf Z}-{\bf U}|^M} \,  \{{\bf Z}, ({\bf Z}-{\bf U})^\dag\} \,   F({\bf Z}) &= F({\bf U}),  & |\underline{y}| &< R.
\end{align*}

\section*{Acknowledgements}
J. Bory Reyes was partially supported by Instituto Polit\'ecnico Nacional in the framework of SIP programs. Alí Guzmán Adán is supported by a BOF-doctoral grant from Ghent University with grant number 01D06014.



\bibliographystyle{abbrv}

\end{document}